\documentclass {article}
\usepackage{amsfonts,amsmath,latexsym}
\usepackage{amstext}
\usepackage {epsf}
\usepackage {epsfig}
\usepackage[scanall]{psfrag}
\usepackage {graphicx}
\usepackage[section]{placeins}
\usepackage {appendix}

\newcommand {\R}{\mathbb{R}}

\newcommand {\h}{\mathbb{H}}
\newcommand {\del}{\partial}

\newcommand {\Span}{\operatorname{span}}

\newcommand {\calV}{\mathcal{V}}
\newcommand {\calL}{\mathcal{L}}
\newcommand {\calW}{\mathcal{W}}

\newtheorem {thm} {Theorem}
\newtheorem {prop} [thm] {Proposition}
\newtheorem {lemma} [thm] {Lemma}
\newtheorem {cor} [thm] {Corollary}
\newtheorem {defn} {Definition}
\newtheorem {rmk} {Remark}

\begin {document}

\title {On the Nondegeneracy of Constant Mean Curvature Surfaces}
\author {Nick Korevaar \\ University of Utah \\
   korevaar@math.utah.edu \and  Rob
   Kusner\footnote{Partially supported by NSF
   grants DMS-0076085 at GANG/UMass and DMS-9810361 at MSRI, and by a 
FUNCAP grant in Fortaleza, Brasil} \\ University of
   Massachusetts, Amherst \\ kusner@math.umass.edu \and Jesse
   Ratzkin\footnote{Partially supported by an NSF VIGRE grant at Utah} \\
University of Connecticut \\ ratzkin@math.uconn.edu }
\maketitle

\begin {abstract}
We prove that many complete, noncompact, constant mean curvature (CMC)
surfaces $f:\Sigma \rightarrow \R^3$ are nondegenerate; that is, the
Jacobi operator $\Delta_f + |A_f|^2$ has no $L^2$ kernel. In fact, if
$\Sigma$ has genus zero with $k$ ends, and if $f(\Sigma)$ is embedded
(or Alexandrov immersed) in a half-space, then we find an explicit upper 
bound for the dimension of the $L^2$ kernel in terms of the number of 
non-cylindrical ends. Our main tool is a conjugation
operation on Jacobi fields which linearizes the conjugate cousin
construction. Consequences include partial regularity for CMC moduli
space, a larger class of CMC surfaces to use in gluing constructions,
and a surprising characterization of CMC surfaces via spinning spheres.
\end {abstract}

\section {Introduction}

Constant mean curvature surfaces in $\R^3$ are equilibria for the area
functional, subject to an enclosed-volume constraint.  The mean
curvature is nonzero when the constraint is in effect, so we can scale
and orient the surfaces to make their mean curvature $1$, a condition
we abbreviate by CMC.  Over the past two decades a great deal of
progress has been made on understanding complete CMC surfaces and
their moduli spaces; however many interesting open problems remain.
One of the most important questions concerns the possibility
of decaying Jacobi fields on complete CMC surfaces, that is, the
Morse-theoretic degeneracy of these equilibria.  The main result
of this paper is to rule out such Jacobi fields on a large class of
complete CMC surfaces.

For a given immersed surface $f: \Sigma \rightarrow \R^3$, 
its mean curvature $H_f$ is determined by the quasilinear elliptic
equation 
$$\Delta_f f = 2H_f\nu_f, $$
where $\nu = \nu_f$ is the (mean curvature, or inner) unit normal to
$f$ and $\Delta_f$ is the 
Laplace-Beltrami operator. The surface $f(\Sigma)$ is CMC if $H_f \equiv 1$. 
The oldest examples of CMC surfaces
are the sphere of radius $1$ and cylinder of radius
$1/2$. Interpolating between these two examples are the Delaunay
unduloids, which are rotationally symmetric and periodic. A
Delaunay unduloid is determined (up to rigid motion) by its
necksize $n$, which is the length of the smallest closed geodesic on
the surface. A necksize of $n = \pi$ corresponds to a cylinder of
radius $1/2$, and as $n \rightarrow 0$ one obtains the singular limit
of a chain of mutually tangent unit spheres. The ODE determining the
Delaunay surfaces still has global solutions when the necksize parameter is
any negative number; in this case the resulting Delaunay nodoids are not
embedded. 

In the present paper we will study CMC 
surfaces in $\R^3$ which are {\it Alexandrov-immersed}. A proper immersion 
$f: \Sigma \rightarrow \R^3$
is an Alexandrov immersion if one can write $\Sigma = \del M$,
where $M$ is a three-manifold into which the mean curvature normal
$\nu$ points, and $f$ extends to a proper immersion of
$M$ into $\R^3$. In the finite topology CMC setting, $M$ is
necessarily a handlebody with a solid cylinder attached for each
end. For example, the Delaunay unduloids are 
Alexandrov-immersed (in fact, embedded), but the 
Delaunay nodoids are not. 

In the remainder of this paper, all of the CMC surfaces are assumed to
be complete, Alexandrov immersions of finite topology, or subsets of
such surfaces.  

It is a theorem of Alexandrov \cite{A} that the only compact, connected, 
Alexandrov-immersed, CMC surfaces are unit spheres. Here we are primarily
interested in noncompact CMC surfaces. Korevaar, Kusner and Solomon
\cite{KKS} proved that each end of such a CMC surface is exponentially
asymptotic to a Delaunay unduloid, that two-ended CMC surfaces are
unduloids, and that three-ended CMC surfaces have  a plane of
reflection symmetry. In fact, all {\em triunduloids} (three-ended,
genus zero CMC surfaces) were constructed and classified by 
Gro\ss e-Brauckmann, Kusner and Sullivan \cite{GKS1}, as were all
{\em coplanar $k$--unduloids} 
($k$--ended, genus zero CMC surfaces whose asymptotic axes all lie in
a plane \cite {GKS2}). These authors define a classifying map
assigning each coplanar $k$--unduloid an 
immersed polygonal disc with $k$ geodesic edges in $S^2$, whose 
edge-lengths are the asymptotic necksizes of the corresponding
$k$--unduloid. 

The classifying map of \cite{GKS1,GKS2} is a homeomorphism, and
gives information about 
the topological structure of moduli space of coplanar
$k$--unduloids. To obtain information about the smooth structure of
moduli space, one needs to understand the 
linearization of the mean curvature operator, which is the Jacobi
operator 
$$
\calL_f = \Delta_f + |A_f|^2,
$$ 
where $|A_f|$ is the length of the second fundamental form of $f$. 
%(The mean curvature operator expresses the first variation of 
%area; the Jacobi operator is its second variation.) 
Solutions to the 
Jacobi equation $\calL_f u = 0$ are called {\em Jacobi fields}, and
correspond to normal variations of the CMC surface $f(\Sigma)$ which
preserve the mean curvature to first order. More precisely, if $u$ is
a Jacobi field, then the one-parameter family of immersions 
$f(t) = f + tu\nu$ has mean curvature $H(t) = 1 + O(t^2)$.
Thus one can think of Jacobi fields as tangent vectors to the moduli
space of constant mean curvature surfaces. 

\begin {defn}
A surface $f:\Sigma \rightarrow \R^3$ is nondegenerate
if the only solution $u \in L^2$ to $\calL_f u = 0$ is the zero function. 
\end {defn}

Near a nondegenerate CMC surface $f(\Sigma)$, a theorem of Kusner,
Mazzeo and Pollack \cite {KMP} shows that the moduli space of CMC
surfaces is a real-analytic manifold with coordinates derived from the
asymptotic 
data of \cite {KKS} (that is, the axes, necksizes, and neckphases of
the unduloid asymptotes). In general the CMC moduli space is a 
real-analytic variety. Indeed, on a degenerate
CMC surface, there would be a nonzero $L^2$ Jacobi field $u$, which
(by \cite{KMP}) decays exponentially on all ends. The presence of such
a Jacobi field means there exists a one-parameter family of surfaces
$f(t)$ with the same asymptotic data and with mean curvature $1 +
O(t^2)$, indicating a possible singularity in the CMC moduli space. 
Thus, proving nondegeneracy eliminates the potential for such
singularities. 

Another application of nondegeneracy is to gluing constructions, where
the pieces to be glued are nondegenerate minimal or CMC surfaces.
Here one looks for CMC surfaces $\hat f(\Sigma)$ which are near a 
given surface $f(\Sigma)$, which may or may not have constant mean 
curvature. Writing the nearby surface $\hat f(\Sigma)$ as a normal
graph over $f(\Sigma)$, with graphing function $u$, we see that
$H_{\hat f} = H_f + \calL_f(u) + O(u^2)$. We take the mean curvature 
of $f$ to be $1 - \phi$, where $\phi$ is a small, but otherwise 
unrestricted, function. Thus, ignoring the higher order terms, finding
the nearby CMC surface $\hat f$ is equivalent to solving the 
linear differential equation $\calL_f(u) = \phi$, with 
growth and/or boundary conditions on $u$. The nondegeneracy of $f$ allows
one to solve this, while controlling $u$ in an appropriate norm. An
early example of this method is Smale's bridge principle \cite{Sm},
which produces a new nondegenerate minimal surface as the
boundary-connected sum of two nondegenerate minimal surfaces along a
thin bridge.
%was the the first gluing result which required nondegenerate
%summands. In this constuction, one joins two minimal (respectively,
%CMC) submanifolds of $\R^N$ using a small bridge, and then perturbs
%this new manifold, which has mean curvature close to $0$ (respectively
%$1$) to obtain a new minimal (CMC) submanifold.  In this case, the
%summands are compact manifolds with boundary, so the nondegeneracy
%condition is slightly different. 
(Kapouleas
\cite{Kap} had earlier constructed many new examples of complete,
noncompact CMC surfaces in $\R^3$ by gluing Delaunay ends and spheres
together.  Although he did not explicitly mention nondegeneracy, his
construction used a balancing condition to overcome the
translation-degeneracy of the sphere.)  Recently, a more flexible
gluing technique \cite{MP, MPP1, MPP2} has been used to explore the moduli
space theory of CMC surfaces: it involves solving several boundary
value problems and then matching Cauchy data across interfaces, but
again the gluing pieces must be suitably nondegenerate.

%RBK - a better transition?
Our main theorem bounds the dimension of the
space of $L^2$ Jacobi fields on a large class of CMC surfaces, and 
as a corollary shows that almost all triunduloids 
are nondegenerate. 

\begin {thm} \label {main-thm}
Let $f:\Sigma \rightarrow \R^3$ be a coplanar $k$--unduloid. 
Then the space of $L^2$ Jacobi 
fields on $f(\Sigma)$ is at most $(k-2)$--dimensional. Moreover, if 
the span of the vertices of the classifying geodesic polygon in $S^2$ 
is $\R^3$, then the space of $L^2$ Jacobi fields on
$f(\Sigma)$ is at most $(k-3)$--dimensional. 
\end {thm}
To precisely state the corollary, 
recall (\cite{GKS1} and our earlier discussion) that a triunduloid
uniquely determines a spherical triangle whose edge-lengths are 
the asymptotic necksizes $n_1,n_2,n_3$. The spherical triangle
inequalities imply $n_1 + n_2 + n_3 \leq 2\pi$ and $n_i + n_j \geq
n_k$. When these inequalities are strict, the vertices of the
classifying triangle span $\R^3$, and so our main theorem asserts that
the space of $L^2$ Jacobi fields is $\{ 0\}$. 

\begin {cor} \label {nondegen-thm}
Let $f:\Sigma \rightarrow \R^3$ be a triunduloid. Then $f$ is
nondegenerate if its necksizes satisfy the strict spherical triangle
inequalities.  
\end {cor}
When a coplanar $k$--unduloid has cylindrical ends, Theorem
\ref{cyl-dim-bound} improves the 
dimension bound in Theorem \ref{main-thm}, and shows many of these CMC
surfaces are also nondegenerate (see Section \ref{questions}). 

The main tool we develop is a conjugate Jacobi field construction,
which converts Neumann fields to Dirichlet fields. This conjugate
variation field arises from linearizing the conjugate cousin
construction of \cite{GKS1}. Although our result is similar in spirit
to the nondegeneracy result of Cos\'in and Ros \cite {CR}, their
method relies on the conformality of the Gauss map, and on the
presence of a homothety Jacobi field.  Both these properties are 
special to the minimal case and do not generalize to our CMC
situation. This difference in the geometry forced us to develop a new
approach.  Not only does it provide another proof of their
nondegeneracy result, but also new insight into the classifying map
for triunduloids and, more generally, for coplanar $k$--unduloids (see
\cite {GKS1, GKS2}). Our conjugate Jacobi field construction also
yields a simple, synthetic characterization of constant mean curvature
in terms of a spinning sphere with speed $2$ along the surface (see
Section \ref{moduli}).

The paper is organized as follows. Notation and preliminary
computations appear in Section \ref{notation}. The conjugate Jacobi
field construction is in Section \ref{conj-var-sect}. In Section
\ref{moduli} we develop the spinning sphere characterization for CMC
surfaces, and the interpretation of the 
classifying map for coplanar $k$--unduloids. 
The proofs of Theorem \ref{main-thm} and Corollary \ref{nondegen-thm}
are in Section \ref{nondegen-sec}. Finally, we discuss some
extensions and applications, as well as  
pose some related open questions, in Section \ref{questions}. 

As with any mathematical problem which has been outstanding for so
long, the present paper has benefited from fruitful 
discussion with many people. In particular, we wish to thank
John Sullivan, Karsten Gro\ss e-Brauckmann and Frank Pacard for reading 
earlier drafts of this paper, and for their helpful suggestions. 

\section {Notation and conventions} \label {notation} 

On a simply connected domain of a CMC (or minimal) surface, we find it
convenient to use {\em conformal curvature coordinates}.
These are coordinates
$(x,y) = (x_1,x_2)$ on a domain $\Omega \subset \R^2$, so that the
mapping $f : \Omega \rightarrow \R^3$ which parameterizes the surface
satisfies  
$$
g_{11}:=\langle f_x,f_x \rangle = \rho^2 = \langle f_y,f_y \rangle=:g_{22}, 
\qquad g_{12}:=\langle f_x,f_y \rangle = 0,
$$
and the (inner) unit normal $\nu$ to the surface satisfies
$$
h_{11}:=\langle \nu,f_{xx} \rangle = 
-\langle \nu_x,f_x \rangle = \rho^2 \kappa_1,
\quad h_{22}:=\langle \nu,f_{yy} \rangle = -\langle \nu_y,f_y \rangle=
\rho^2 \kappa_2, 
\quad h_{12}:=\langle \nu,f_{xy} \rangle = 0.
$$ In other words, choosing conformal curvature coordinates amounts to
simultaneously diagonalizing the first and second fundamental forms,
$g$ and $h$. In these coordinates, the shape operator $A = g^{-1} h$
is diagonal with the principal curvatures $\kappa_1, \kappa_2$ as its
entries. Equivalently, the $x$ and $y$  coordinate lines are principal
curves. Notice that $H = (\kappa_1 + \kappa_2)/2$ is half the trace of
$A$. In what follows it will be useful to define  $\kappa:=
\kappa_2 - \kappa_1$, and to adopt the convention $\kappa_2 > 
\kappa_1$ (away from umbilic points). 
It also will be convenient to decompose the shape operator as $A = B + C$,
where $C = HI$ is the trace part and $B = A-HI$ is trace-free. Thus, in
conformal curvature coordinates, $A$ and $B$ have matrices
$$ A = 
\left [ \begin {array}{cc} \kappa_1 & 0 \\ 0 & \kappa_2 \end {array}
  \right ], 
\qquad B = \left [ \begin {array}{cc} (\kappa_1 - \kappa_2)/2 & 0
    \\ 0 & (\kappa_2 - \kappa_1)/2 \end {array} \right ]
= \left [ \begin {array}{cc} -\kappa/2 & 0
    \\ 0 & \kappa/2 \end {array} \right ].
$$ 

The existence of conformal curvature coordinates (away from
umbilics) on a CMC surface can be seen using the {\em Hopf differential}, a
holomorphic quadratic differential associated with $B$ (see
\cite{Ho}).  More precisely, suppose we 
have any conformal coordinates $(u,v)$ on the surface, and consider 
the complex coordinate $w=u+iv$. The Codazzi equation 
implies the complex-valued function
$$
\phi: = (h_{11}-h_{22})/2 +ih_{12}
$$ 
is holomorphic with respect to $w$ if
and only if $H$ is constant. Under conformal changes of coordinates,
the holomorphic quadratic differential
$$ \Phi:= \phi(w) dw^2 $$  
is invariant.  This $\Phi$ is the Hopf differential of our CMC surface.  

\begin {lemma}  \label {hopf-const}
If $\Omega$ is simply connected and $f : \Omega \rightarrow \R^3$ is a
conformal immersion of a CMC surface without umbilics, then
there exists a conformal change of coordinates so that $f$ is an
immersion with conformal curvature coordinates, and so that
$\kappa>0$. Moreover, in any conformal curvature coordinates, $\kappa
\rho^2$ is a constant. 
\end {lemma}

{\bf Proof}: Observe that umbilic points of $f$ are precisely the
zeroes of $\Phi = \phi(w) dw^2$. Because 
$\Omega$ is simply connected and $f(\Omega)$ has no umbilics, we can
pick a branch of $\sqrt{\phi(w)}$. Make a conformal change of  
coordinates $z = z(w) = x + iy$ by integrating the one-form 
\begin {equation} \label {hopf-const-eqn}
dz := i\sqrt{\Phi} = i\sqrt{\phi(w)} dw. \end {equation}
Then in the $z$ coordinates, $\Phi = -dz^2$. 
This means $h_{12} \equiv 0$, and so $f(w(z))$ is an immersion in conformal
curvature coordinates. Also, 
$h_{11}- h_{22} = -2$ implies $\kappa \rho^2 =2$, and so
$\kappa> 0$. 

Moreover, for any conformal curvature coordinates, $h_{12} \equiv 0$,
so $-2\phi = \kappa \rho^2$ is a real-valued holomorphic function, and
hence constant.  
\hfill $\square$ 

\hskip .1in

We now proceed with some preliminary computations using conformal
curvature coordinates. These are
elementary, but we include them for the convenience of the
reader.  Using the flat Laplacian, $\Delta_0 = \del_x^2 + \del_y^2$,
the CMC equation is 
$$\rho^2 \Delta_f f = \Delta_0 f = 2f_x \times f_y = 2\rho^2 \nu,$$  
and the Jacobi equation reads
\begin {equation} \label {jacobi-eqn}
\rho^2 \calL_f u = \Delta_0 u + \rho^2(\kappa_1^2 + \kappa_2^2) u = 0.
\end {equation}

Unlike the previous lemma, the next two do not require $f$ to have
constant mean curvature. However, they do require that $f:\Omega
\rightarrow \R^3$ is an immersion in conformal curvature coordinates. 

\begin {lemma} \label {2nd-deriv-f}
If $f:\Omega \rightarrow \R^3$ is an immersion in conformal curvature
coordinates, with unit
normal $\nu$ and conformal factor $\rho=|f_x|=|f_y|$, then we have
$$f_{xx} = \frac{\rho_x}{\rho} f_x - \frac{\rho_y}{\rho} f_y +
  \kappa_1 \rho^2 \nu, \qquad f_{yy} = -\frac{\rho_x}{\rho} f_x +
  \frac{\rho_y}{\rho} f_y + \kappa_2 \rho^2 \nu.$$
\end {lemma}

{\bf Proof}: The frame $(f_x,f_y,\nu)$ is orthogonal, so 
$$f_{xx} = \rho^{-2} \langle f_{xx}, f_x \rangle f_x + \rho^{-2} \langle
f_{xx}, f_y \rangle f_y + \langle f_{xx}, \nu \rangle \nu, \qquad 
f_{yy} = \rho^{-2} \langle f_{yy}, f_x \rangle f_x + \rho^{-2} \langle
f_{yy}, f_y \rangle f_y + \langle f_{yy}, \nu \rangle \nu.$$
One can then complete the proof by differentiating the equations
$$\langle f_x,f_x \rangle = \rho^2 = \langle f_y, f_y \rangle, \quad
\langle f_x, f_y \rangle = \langle f_x, \nu \rangle = \langle f_y, \nu
\rangle = 0.$$
\hfill $\square$ 

\begin {lemma} \label {var-cmplx-struct}
If $f: \Omega \rightarrow \R^3$ is an immersion in conformal curvature
coordinates and $u \in C^2(\Omega)$, then one can write
the complex structure of the surface $f(t) = f+tu\nu + O(t^2)$ 
as $$J(t) = J_0  + tJ_1 + O(t^2)$$ where
$$J_0 = \left [ \begin {array}{cc} 0 & -1 \\ 1 & 0 \end {array} \right
], \qquad J_1 = \left [ \begin {array}{cc} 0 & u\kappa
    \\ u\kappa & 0\end {array} \right ].$$
Thus the coordinate-free expression for $J_1$ is the product $2u J_0
B$, where $B$ is the trace-free shape operator of $f$.
\end {lemma}

{\bf Proof}: In any oriented local coordinates, 
$$J = \frac{1}{\sqrt {\det(g)}} \left [ \begin {array}{cc} -g_{12} & -
    g_{22} \\ g_{11} & g_{12} \end {array} \right ].$$
Using conformal curvature coordinates at $t=0$, we compute the metric
    at $t$:
\begin {eqnarray*}
g_{11} & = & \langle f_x(t) ,f_x(t) \rangle = \rho^2 (1 - 2tu\kappa_1)
+ O(t^2)\\  
g_{22} & = & \langle f_y(t) ,f_y(t) \rangle = \rho^2 (1 - 2tu\kappa_2)
+ O(t^2)\\  
g_{12} & = & \langle f_x(t) ,f_y(t) \rangle = O(t^2). \end {eqnarray*}
Thus 
\begin {eqnarray*}
J(t) & = & \frac{1}{\rho^2 \sqrt{1 - 2(\kappa_1 + \kappa_2) tu}} \left
    [ \begin {array}{cc} 0 
    & -\rho^2 (1 - 2tu\kappa_2) \\ \rho^2(1 - 2tu\kappa_1) & 0 \end
    {array} \right ] + O(t^2) \\
& = & ( 1 + tu(\kappa_1 + \kappa_2)) \left [ \begin {array}{cc} 0 & -1 +
    2tu\kappa_2 \\ 1 - 2tu\kappa_1 & 0 \end {array} \right ] + O(t^2), 
\end {eqnarray*}
which yields the desired expansion. 
\hfill $\square$
\vskip .1in

Lawson \cite{L} pioneered the conjugate cousin relation between CMC
surfaces and minimal surfaces in $S^3$.
The first order conjugate cousin construction was
initiated by Karcher \cite {K} and developed in \cite{G,GKS1}. It
uses the realization of $S^3 \subset \R^4 = \h$ as the unit
quaternions, and of $\R^3 = \Im \h$ (the
imaginary quaternions) as the Lie algebra of $S^3$, or
as the tangent space $T_1 S^3$. For imaginary quaternions 
$p, q \in \R^3$, we can write their product as 
\begin {equation} \label {cross-product}
p q = -\langle p,q \rangle + p \times q. \end {equation}
In particular, orthogonal imaginary quaternions anti-commute.
Thus we can also write the CMC condition $H_f \equiv 1$ as 
\begin {equation} \label {cmc-eqn}
\Delta_0 f = 2f_x f_y = 2\rho^2 \nu. \end {equation}

Let $\Omega \subset \R^3$ be a simply connected domain. Theorem 1.1 of
\cite{GKS1} shows that conjugate
cousins $f : \Omega \rightarrow \R^3$ and $\tilde f : \Omega
\rightarrow S^3$ satisy the first order system of partial differential
equations 
\begin {equation} \label {conj-cousin} 
d \tilde f = \tilde f df \circ J_0, \end {equation}
where $J_0$ is the standard complex structure on $\R^2$ and the product
is the quaternion product. The integrability condition for $\tilde f$
reduces to the CMC equation for $f$, and in this case the resulting
surface $\tilde f(\Omega) \subset S^3$ is minimal. Conversely, given a
minimal immersion $\tilde f$, one can
consider $f$ as the unknown in the system (\ref{conj-cousin}). Then
the integrability condition for $f$ is the minimality of $\tilde f$,
and the resulting surface $f(\Omega) \subset \R^3$ is CMC. Moreover,
the immersions 
$f$ and $\tilde f$ are uniquely determined up to translation in
$\R^3$ and left translation in $S^3$, respectively. One can also see
from equation (\ref{conj-cousin}) that $f$ and $\tilde f$ are
isometric. 

\begin {lemma} \label{same-Jacobi-op}
The Jacobi operators for $f$ and $\tilde f$ coincide, and so we can
identify Jacobi fields on the two surfaces.
\end {lemma}

{\bf Proof}: In general, the Jacobi operator for a two-sided (CMC or
minimal) surface
with normal $\nu$ in a manifold with Ricci 
curvature $\operatorname{Ric}$ is 
$$\calL = \Delta + |A|^2 + \operatorname{Ric}(\nu, \nu).$$
Since the Ricci curvature of $\R^3$ or $S^3$ is $0$ or $2$,
respectively, for $f$ and its cousin $\tilde f$ we have 
$$\calL_f = \Delta_f + |A_f|^2, \qquad \calL_{\tilde f} =
\Delta_{\tilde f} + |A_{\tilde f}|^2 + 2. $$
The two surfaces are isometric, so $\Delta_f = \Delta_{\tilde
  f}$. Moreover, we have (see Proposition 1.2 of \cite {GKS1}) 
$$\tilde \kappa_1 = \kappa_1 - 1, \qquad \tilde \kappa_2 = \kappa_2
-1.$$
Thus 
$$|A_{\tilde f}|^2 = \tilde \kappa_1^2 + \tilde \kappa_2^2 =
(\kappa_1-1)^2 + (\kappa_2 - 1)^2 = \kappa_1^2 + \kappa_2^2 -
2(\kappa_1 + \kappa_2) + 2 =
|A_f|^2 -2.$$
\hfill $\square$

\section {Existence of the conjugate cousin variation field}
\label {conj-var-sect}

In this section we construct a conjugate cousin variation field
$\tilde \epsilon$ on $\tilde f$ from a normal variation field
$u\nu$ on $f$. The idea behind this construction is to linearize the
conjugate cousin equation (\ref{conj-cousin}). 

We begin with a CMC immersion 
$f: \Omega \rightarrow \R^3$ of a simply connected domain 
and a solution 
$u: \Omega \rightarrow \R$
to the Jacobi equation (\ref{jacobi-eqn}).
In general, $u\nu$ is not the initial velocity of a one-parameter
family of CMC surfaces 
$$f(t) = f + tu\nu + O(t^2).$$
Although one can always find such a family on a sufficiently small
subdomain, 
the families will not coincide on the overlaps of these
subdomains. However, when there does exist such a one-parameter family of
CMC surfaces, then one can define a conjugate cousin family 
$$\tilde f(t) = \tilde f + t\tilde \epsilon + O(t^2)$$
by integrating equation (\ref{conj-cousin}). 
In this case, the two families are related by the system  
\begin {equation} \label {conj-family}
d \tilde f(t) = \tilde f (t) df(t) \circ J(t), \end {equation}
where $J(t)$ is the complex structure on $f(t)$. Surprisingly,
if the domain $\Omega$ is simply connected, then an initial
velocity field $\tilde \epsilon$ can be defined globally, even though
this may not be possible for the conjugate cousin family itself. 

\begin {prop} \label{conj-var-exist}
Let $p$ be a point in a simply connected domain $\Omega$. Let $f$
 and $\tilde f$ be conjugate cousins satisfying equation
 (\ref{conj-cousin}). Then for any Jacobi field $u$ on
 $\Omega$, and any 
choice of initial velocity $\tilde\epsilon(p)$,  there exists
a unique global variation field $\tilde \epsilon$ on $\tilde f(\Omega)$
which is locally associated to $u$ in the manner  described
above.  The field $\tilde \epsilon$ satisfies the first order system
of linear partial differential equations
\begin {equation} \label {conj-variation}
d\tilde \epsilon = \tilde f df \circ J_1 + \tilde f d(u\nu) \circ J_0
+ \tilde \epsilon df \circ J_0. \end {equation}
\end{prop}

\begin{rmk} The new variation field $\tilde \epsilon$ need not be a 
 normal field along $\tilde f$. \end {rmk}

{\bf Proof:} We first sketch an abstract proof of the
proposition, before giving a purely computational one.
Small patches
of a CMC surface are graphical and therefore strictly stable.
Thus one can always use the implicit function
theorem to solve a family of Dirichlet problems for the normal
variation CMC equation, with boundary data $f(t) = f + tu\nu$.
This yields a one-parameter family of CMC patches $f(t)$ with $t$
in a neighborhood of $0$, and with 
initial velocity $u\nu$ on such a small patch.
From these CMC patches, solve equation (\ref{conj-family}) for a
family $\tilde f(t)$ of  minimal
surface patches in $S^3$, uniquely determined for each $t$ once one 
specifies a basepoint $\tilde\gamma(t)=\tilde f(t)(p)$.  These conjugate
cousin surfaces have an initial velocity field $\tilde \epsilon$. 
Note, $\tilde\epsilon(p)=\tilde\gamma'(0)$
can be adjusted at will. Once we show that the fields $\tilde\epsilon$
all satisfy the first order system (\ref{conj-variation}) we deduce
 not only local existence for the initial value
problem (as just described), but also
uniqueness, since equation (\ref{conj-cousin}) reduces to a first order
system of differential equations along any curve. 
Global existence and uniqueness then follow because $\Omega$ is
simply connected.

To derive our governing system (\ref{conj-variation})
 we expand the conjugate family equation (\ref{conj-family}) (using
quaternionic multiplication throughout):
\begin {eqnarray*} 
d \tilde f + t d \tilde \epsilon + O(t^2) 
& = & 
d \tilde f(t) 
 = 
\tilde f(t) df(t) \circ J(t) \\ 
& = & 
(\tilde f + t\tilde \epsilon + O(t^2)) (df + td(u\nu) + O(t^2))
\circ (J_0 + tJ_1 + O(t^2)) \end {eqnarray*}
Equating the $O(1)$ terms in this expansion gives the cousin equation
(\ref{conj-cousin}).
  Equating the $O(t)$ terms yields equation (\ref{conj-variation}),
  completing our sketch of the abstract proof.

A direct and instructive proof of Proposition \ref{conj-var-exist}
is to show that the first order system of partial differential
equations (\ref{conj-variation}) satisfies the Frobenius
integrability conditions, namely that the formal mixed partial derivatives
are equal.  Existence and uniqueness for the initial value problem then
follows directly from the Frobenius theorem and the fact that $\Omega$
is simply connected.  Verifying the mixed-partials condition
 amounts to showing
that the formal computation of
$d(d\tilde \epsilon)$ yields $0$. Differentiating and expanding
equation (\ref{conj-variation}),
we get eight terms: 
\begin {eqnarray} \label {mixed-partials}
d(d\tilde \epsilon) & = & d(\tilde f df \circ J_1) + d(\tilde f
d(u\nu) \circ J_0) + 
d\tilde \epsilon \wedge df \circ J_0 + \tilde \epsilon d(df \circ J_0)
\\ \nonumber 
& = & \tilde f df\circ J_0 \wedge df \circ J_1 + \tilde f d(df \circ
J_1) + \tilde f df \circ J_0 \wedge d(u\nu) \circ J_0 + \tilde f
d(d(u\nu) \circ J_0) \\ \nonumber
& & + \tilde f df \circ J_1 \wedge df \circ J_0+
\tilde f d(u\nu) \circ J_0 \wedge df \circ J_0 + \tilde \epsilon df
\circ J_0 \wedge df \circ J_0 + \tilde \epsilon d(df \circ J_0).
\end {eqnarray}

It is easiest to analyze equation (\ref{mixed-partials}) term by term.
We use conformal curvature coordinates to compute coordinate-free
identities. Since umbilic points are isolated (we are not considering
subdomains of spheres), continuity implies these identities 
hold everywhere. 
All terms are multiples of the area form $da = \rho^2 dx\wedge dy$, and
two of the terms vanish:

\begin {lemma}
$$df \circ J_1 \wedge df \circ J_0 = 0 = df \circ J_0 \wedge df \circ
  J_1.$$
\end {lemma}

{\bf Proof}: We compute $df\circ J_1 \wedge df \circ J_0$: 
\begin {eqnarray*}
df \circ J_1 \wedge df \circ J_0 & = & (u\kappa f_y dx + u\kappa f_x)
\wedge (f_y dx - f_x dy) = u\kappa(-f_y f_x dx \wedge dy + f_x f_y dy
\wedge dx) \\ 
& = & u\kappa (f_x f_y - f_x f_y) dx \wedge dy = 0. 
\end {eqnarray*}
Here $f_x$ and $f_y$ are orthogonal, so they
anti-commute. We also have 
$$df \circ J_0 \wedge df \circ J_1 = - df \circ J_1 \wedge df \circ
J_0 = 0.$$
\hfill $\square$

Using equation (\ref{cmc-eqn}), the next lemma implies that two more
terms sum to zero:
\begin {lemma}
$$d(df \circ J_0) = - \Delta_0 f dx \wedge dy = -2\rho^2 \nu dx
  \wedge dy = -2\nu da, \qquad df \circ J_0 \wedge df \circ J_0 =
  2\rho^2 \nu dx 
  \wedge dy = 2 \nu da.$$
\end {lemma}

{\bf Proof}: First we compute
$$ d(df \circ J_0) = d(f_y dx - f_x dy) = f_{yy} dy \wedge dx -
f_{xx} dx \wedge dy = - \Delta_0 f dx \wedge dy = -2\rho^2 \nu dx
\wedge dy.$$
Similarly, 
$$df \circ J_0 \wedge df \circ J_0 = (f_y dx - f_x dy) \wedge (f_y dx
- f_x dy) = -f_y f_x dx \wedge dy - f_x f_y dy \wedge dx = 2\rho^2 \nu
dx \wedge dy.$$ 
\hfill $\square$

The remaining terms involve the decomposition of the shape operator $A$ 
into trace-free
and trace parts, $B = A - C$ and $C= HI = I$, respectively.  In fact,
note that $A = B + C$ is an orthogonal decomposition 
in the space of symmetric linear maps, so that, by the Pythagorean theorem,
$$
|A|^2 = |B|^2 + |C|^2.
$$

\begin {lemma}
$$
d(df \circ J_1) = -2[df (B\nabla u) + |B|^2 u\nu] da. 
$$
\end {lemma}

{\bf Proof}: We compute, using Lemmas \ref{2nd-deriv-f} and
\ref{hopf-const}: 
\begin {eqnarray*} 
d(df \circ J_1) & = & d(u\kappa f_y dx + u\kappa f_x dy) \\ 
& = & [u_x \kappa f_x + u \kappa_x f_x + u\kappa f_{xx}]dx \wedge dy +
[-u_y \kappa f_y - u\kappa_y f_y - u\kappa f_{yy}] dy \wedge dx \\ 
& = & [\kappa(u_x f_x - u_y f_y) + u(\kappa_x f_x - \kappa_y f_y) +
  u\kappa (f_{xx} - f_{yy})] dx \wedge dy \\ 
& = & [\kappa(u_x f_x - u_y f_y) + u(\kappa_x f_f - \kappa_y f_y) +
  u\kappa (2\rho^{-1} \rho_x f_x - 2\rho^{-1} \rho_y f_y- \kappa
  \rho^2 \nu) ] dx \wedge dy \\
& = & [\kappa(u_x f_x - u_y f_y) + u(\kappa_x f_x + 2\kappa \rho^{-1}
  \rho_x f_x - \kappa_y f_y- 2\kappa \rho^{-1} \rho_y f_y - \kappa^2
  \rho^2 \nu )]dx \wedge dy \\ 
& = & [\kappa(u_x f_x - u_y f_y) + u(2\rho^{-2}\del_x(\kappa \rho^2)
  f_x - 2\rho^{-2} \del_y(\kappa \rho^2) f_y - \kappa^2 \rho^2 \nu)]
dx \wedge dy \\ 
& = & [\kappa(u_x f_x - u_y f_y) - \kappa^2 \rho^2 u\nu] dx \wedge
dy  =  -2[df(B\nabla u) + |B|^2 u\nu] da. 
\end {eqnarray*}
\hfill $\square$

The next term we have is:
\begin {lemma}
$$
d(d(u\nu) \circ J_0) =  2[ df(A\nabla u) + |A|^2 u\nu] da.
$$
\end {lemma} 
{\bf Proof}: We compute, using the Jacobi equation:
\begin {eqnarray*} 
d(d(u\nu) \circ J_0) & = & d((u\nu)_y dx - (u\nu)_x) dy) = - \Delta_0
(u\nu) dx \wedge dy \\ 
& = & -[u \Delta_0 \nu + (\Delta_0 u)\nu + 2\langle \nabla u, \nabla
  \nu \rangle] dx \wedge dy \\ 
& = &  - [-u\rho^2 |A|^2 \nu -\rho^2 |A|^2 u\nu + 2u_x \nu_x + 2u_y
  \nu_y] dx \wedge dy \\ 
& = & (2\kappa_1 u_x f_x + 2\kappa_2 u_y f_y + 2\rho^2 |A|^2 u\nu) dx
\wedge dy  =  2[ df(A\nabla u) + |A|^2 u\nu] da.
\end {eqnarray*} 
\hfill $\square$

The final two terms actually coincide: 

\begin {lemma}
$$d(u\nu) \circ J_0 \wedge df \circ J_0 = 
 -[df(C\nabla u) + |C|^2 u\nu] da = df \circ J_0 \wedge d(u\nu) \circ
  J_0. $$ 
\end {lemma}
{\bf Proof}: Using the conformality relations $\nu f_x = f_y$ and
$\nu f_y = -f_x$, we have
\begin {eqnarray*}
d(u\nu) \circ J_0 \wedge df \circ J_0 & = & ((u_y \nu + u \nu_y)dx -
(u_x \nu + u \nu_x) dy) \wedge (f_y dx - f_x dy) \\
& = & ((u_y \nu -
\kappa_2 u f_y)dx - (u_x \nu - \kappa_1 f_x)dy) \wedge (f_y dx - f_x
dy) \\
& = & (-u_y \nu f_x + \kappa_2 u f_y f_x) dx \wedge dy + (-u_x \nu f_y
+ \kappa_1 u f_x f_y) dy \wedge dx \\ 
& = & (-u_y f_y - \kappa_2 \rho^2 u\nu)
dx \wedge dy + (u_x f_x + \kappa_1 \rho^2 u\nu) dy \wedge dx \\ 
& = & (-u_x f_x - u_y f_y - (\kappa_1 + \kappa_2) \rho^2 u\nu) dx
\wedge dy \\ 
& = & (-u_x f_x - u_y f_y - 2\rho^2 u\nu) dx \wedge dy
 =  -[ df(C\nabla u) + |C|^2 u\nu] da,
\end {eqnarray*}
since CMC implies the trace-part $C = I$ and thus $|C|^2 = 2$.
The other computation is similar.
\hfill $\square$ 
\vskip .1in

Summing the results of the previous lemmas:
$$
d(d\tilde \epsilon) 
= 2\tilde f[ df((A-B-C)\nabla u) + (|A|^2-|B|^2-|C|^2) u\nu] da
= 2\tilde f[0 + 0] = 0. 
$$
This completes the proof of the proposition. 
\hfill $\square$

\section {Homogeneous solutions, spinning spheres, and the classifying
  map via pole solutions}
\label {moduli}

We continue to consider a simply connected CMC surface $f:\Omega
\rightarrow \R^3$ and its conjugate cousin 
surface $\tilde f: \Omega\rightarrow S^3$.
At this point, it is useful to pull the variation field $\tilde
\epsilon$ back to $\R^3=T_1 S^3$. Thus we define 
$$\epsilon := \tilde f^{-1} \tilde \epsilon.$$
By the product rule and equation (\ref{conj-cousin}), we have 
$$
d\tilde \epsilon = d(\tilde f \epsilon) = \tilde f (df \circ J_0)
\epsilon + \tilde f d\epsilon;
$$
however, by equation (\ref{conj-variation}), 
$$
d\tilde \epsilon = \tilde f df \circ J_1 + \tilde f d(u\nu) \circ J_0
+ \tilde \epsilon df \circ J_0 = \tilde f (df \circ J_1 + d(u\nu)\circ
J_0 + \epsilon df \circ J_0).
$$
Equating these two expressions, solving for $d\epsilon$, and
applying equation (\ref{cross-product}), one obtains
\begin {equation} \label {revised-conj-variation}
d\epsilon  = \epsilon (df \circ J_0) - (df \circ J_0) \epsilon + df \circ
J_1 + d(u\nu) \circ J_0 
=  2\epsilon \times df \circ J_0 + df \circ J_1
+ d(u\nu) \circ J_0. \end {equation}

\subsection{Homogeneous solutions and spinning spheres} 
\label{spin-sphere-sec}

Equation (\ref{revised-conj-variation}) is an inhomogeneous first
order differential system for $\epsilon$, where the inhomogeneity $df
\circ J_1 + d(u\nu) \circ J_0 = (-2df \circ Bu + d(u\nu)) \circ J_0$ 
depends
linearly on the Jacobi field $u$. When $u \in L^2$, the asymptotic behavior of
a solution to equation (\ref{revised-conj-variation}) is given by solutions to 
the associated homogeneous ($u \equiv 0$) equation: 
\begin {equation} \label {compass-eqn}
d\epsilon = \epsilon (df \circ J_0) - (df \circ J_0) \epsilon =
  2\epsilon \times df \circ J_0.
\end {equation}

We first study the geometry of solutions to equation (\ref{compass-eqn}). 
Notice that equation (\ref{compass-eqn}) implies $\epsilon$ is
perpendicular to $d\epsilon$, so
$$d(|\epsilon|^2) = 2\langle d\epsilon, \epsilon \rangle = 0,$$
and the solutions  $\epsilon$ to the linear system (\ref{compass-eqn})
have  globally constant length. It follows that one can
use them to define a path-independent parallel transport along 
$f(\Omega)$,
mapping $T_{f(p)}\R^3\rightarrow T_{f(q)}\R^3$ isometrically.
To see this, let $\gamma$ be path from $p$ to $q$ on the simply
connected domain $\Omega$. One recovers $\epsilon(f(q))$ by
integrating the solution to the initial value problem for equation
(\ref{compass-eqn}), with initial value $\epsilon_0 = \epsilon(f(p))$.
Since this parallel transport is path independent,
it defines a flat connection on a principal $SO(3)$-bundle over 
$\Omega$.

There is an interesting physical interpretation of this flat
connection.  Notice that if one integrates equation
(\ref{compass-eqn}) along any curve $\gamma$ then a solution
$\epsilon$ with unit length rotates with constant angular speed $2$, 
with
evolving axis of rotation given by the curve conormal, $df \circ
J_0(\gamma'(s))=\eta(s)$.  This means that the $SO(3)$-frame evolves
as if it were attached to a sphere spinning around the conormal $\eta$
at speed $2$.
%Note that the way to roll a sphere
%along the surface, without twisting or slipping, so that total
%rotation is minimized, is to have the
%sphere rotate about an axis parallel to the contact curve
%conormal. (Precisely, the total rotation is the length of a path in
%$SO(3)$, which we minimize subject to the constraint that the spinning
%sphere is tangent to $f(\Omega)$ as it traverses $f(\gamma)$.) If
%the sphere has radius $1/2$, and the contact point moves at speed $1$,
%then the angular speed of rotation is $2$.
%If the
%surface $f(\Omega)$ has mean curvature $1$, and if the radius--$1/2$
%sphere is on the outside of the surface relative to the inner
%normal $\nu$, then the spinning sphere exactly reproduces our
%flat connection. (One must allow the sphere to
%immerse through the surface as necessary, for example near points with
%a principal curvature less that $-2$; in fact, the sphere
%should really roll with axis tangent to the CMC surface, but that
%equivalent spinning motion would be impossible to carry out
%physically.)
In particular, if the spinning sphere follows a
(contractible)
loop on the
surface, it will return with its initial orientation. This even gives
a surprising property on a round sphere. %The physical
%realization of this mathematical fact would make an interesting
%demonstration. 

In fact, the flatness of this connection is equivalent to $f$ having 
mean curvature $1$. 
\begin {prop} \label{spin-sphere-prop}
Let $f:\Omega \rightarrow \R^3$ be an immersion and consider the
$SO(3)$-connection defined by spinning a sphere at speed $2$ as
described above. Then $f(\Omega)$ has mean curvature $1$ if and only
if this connection is flat. 
\end {prop}

We prove this proposition and further explore the spinning sphere 
connection in 
Appendix \ref{spin-sphere-appen}. 

\subsection{Pole solutions and the classifying map}\label{modulisub} 

The $\epsilon$-fields which solve the homogeneous system (\ref{compass-eqn})
yield a new perspective on the classifying map \cite {GKS1, GKS2}
for coplanar $k$--unduloids. 

\begin {figure}[h]
\begin {center}
\includegraphics[width=4in]{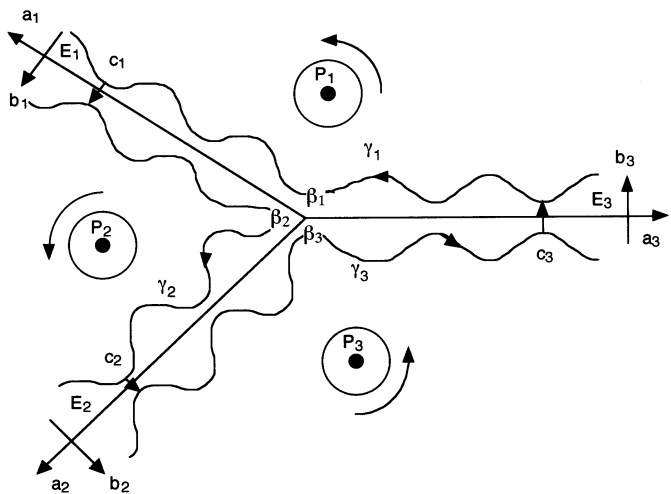}
\caption{Triunduloid configuration from above}
\label {configuration}
\end {center}
\end {figure}

Let $f:\Sigma \rightarrow \R^3$ be a CMC surface with $k$ ends and 
genus zero, which lies in a half-space (necessarily so when 
$k \leq 3$). By \cite {KKS} so a coplanar $k$--unduloid is 
{\em Alexandrov symmetric}: it has a reflection plane of symmetry,
which we normalize to be the $xy$ plane; furthermore, the closures of
each half of $f(\Sigma)$, $f(\Sigma^+) = f(\Sigma) \cap \{z > 0 \}$
and $f(\Sigma^-) = f(\Sigma) \cap \{ z < 0\}$, are graphs over a
(possibly immersed) planar domain. Because $\Sigma$ has genus zero,
$\Sigma^\pm$ are topological discs. The common boundary $\del
f(\Sigma^\pm)$ is the union of $k$ oriented, planar, principal curves
$\gamma_1, \dots, \gamma_k$, where $\gamma_j$ connects the end
$E_{j-1}$ to $E_j$, using the natural cyclic ordering of the ends (see
\cite{GKS2}). The configuration for a triunduloid ($k = 3$) is
indicated in Figure \ref{configuration}.

The evolution of solutions to equation (\ref{compass-eqn}) is easy to
track along 
curves of constant conormal $\eta(s)=df\circ J_0(\gamma'(s))$, 
since the conormal is the rotation axis.
With our convention that the inner normal $\nu = \gamma '(s)\eta(s) 
= \gamma '(s)\times\eta(s) $,
and our choice of curve orientation in Figure \ref{configuration},
we see that the rotation axis along each $\gamma_j$ is the vertical 
vector $\eta=-e_3$, so that the rotation appears counterclockwise
from above, as indicated in the figure.  

\begin {figure}[h]
\begin {center}
\includegraphics[width=2.5in]{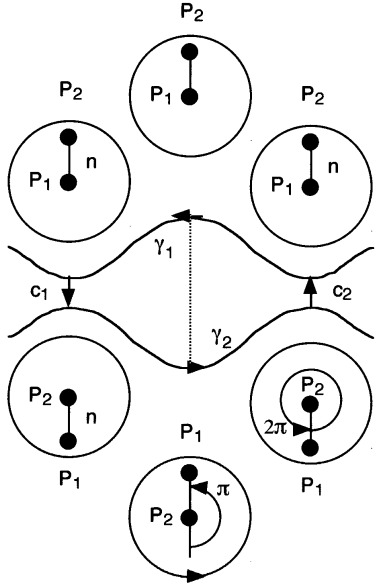}
\caption{Delaunay configuration from above}
\label {delaunay-contour}
\end {center}
\end {figure}

\begin {defn} 
The pole solutions $P_1, \dots, P_k$ to equation (\ref{compass-eqn}) are 
solutions 
to the ODE system on $f(\bar \Sigma^+)$, with the initial value $P_j = e_3$ at 
some point (and hence all points) of $\gamma_j$. 
\end {defn}
Observe that each $P_j$ is a globally defined unit vector field on 
$f(\bar \Sigma^+)$. Moreover, equation (\ref{compass-eqn}) shows that the
angles between any two pole solutions $P_i$ and $P_j$ remains constant. 
Thus, along any curve $\gamma$, all pole solutions evolve by the same 
rotation, and so $P_1, \dots, P_k$ can be viewed as the vertices of a 
geodesic polygon in $S^2$ which is well-defined up to rotation. We begin 
with a lemma about the pole solutions in the case $k=2$, which we will 
also need in Section \ref{neumann-sec}.
\begin{lemma}\label{delaunay} 
Let $f(\Sigma)$ be an unduloid with profile curves $\gamma_1, 
\gamma_2$ as in Figure \ref{delaunay-contour}.
\begin {itemize}
\item  If $f(\Sigma)$ is not a cylinder then $\gamma_1$ and $\gamma_2$ each 
have period $\pi$ when parameterized by arclength (see also Section 1 of 
\cite{GKS1}).  
%\item  The conjugate surface
%$\tilde f(\Sigma^+)$ 
%  periodically covers an embedded minimal annulus in $S^3.$ 
\item If $f(\Sigma)$ is non-cylindrical
then the only solution on $f(\Sigma^+)$ to equation (\ref{compass-eqn}) 
which satisfies
 $\langle \epsilon,\nu\rangle=0$  along
both $\gamma_1,\gamma_2$ is the zero solution.  If $f(\Sigma)$
is cylindrical, then on $f(\Sigma^+)$ the pole solutions $P_1, P_2$ are 
opposites, and
are tangential to $f(\Sigma^+)$, that is $\langle P_j, 
\nu\rangle \equiv 0$.
Each solution $\epsilon$ of equation (\ref{compass-eqn}) satisfying
$\langle \epsilon,\nu\rangle=0$  along both $\gamma_1,\gamma_2$ 
is a multiple of $P_1 = -P_2$.
\end{itemize}
\end{lemma}

{\bf Proof:} 
Starting at the initial point of $c_1$ in Figure \ref{delaunay-contour}, 
follow the pole solutions
around the contour in thie figure, which depicts one
period of an 
unduloid. We see that the pole $P_1$ must 
return to the vertical position after traversing the second neck
$c_2$. This is only possible if $P_1$ has rotated through a total angle
of $2\pi k$ for
some positive integer $k$ as it travels from $c_1$ to $c_2$ along
$\gamma_2$. However, $P_1$ 
rotates with speed $2$ along $\gamma_2$, so the length of the
$\gamma_2$--arc 
must be $k\pi$. In the zero necksize limit, this arc is half a great
circle on a unit sphere, so it has length $\pi$. Thus, by the
continuity of the family of Delaunay unduloids, the period of each
unduloid is $\pi$.  

%The fact that $\tilde f(\Sigma^+)$ covers an embedded minimal surface
%follows from the fact that the profile curves of $f$ have period
%$\pi$. Indeed, these curves have constant conormal, so their cousin
%curves are great circle arcs of length $\pi$. Thus traversing $2$
%periods on the Delaunay unduloid will return one to the beginning of
%the cousin great circles, with the same tangent direction, showing
%that the cousin minimal surface is embedded. 

For the second part of this proposition, suppose $\epsilon \neq 0$
solves equation (\ref{compass-eqn}) and $\langle \epsilon, \nu \rangle
= 0$ along $\gamma_1,\gamma_2$. If $\epsilon$ has a nonzero
horizontal component along the boundary 
curve $\gamma_1$, then as one traverses $\gamma_1$ this component 
rotates with angular speed $2$. Thus the horizontal component of
$\epsilon$ will be perpendicular to the axis of the unduloid at points
distributed with period $\pi/2$. At such points $\langle \epsilon, \nu
\rangle \neq 0$. Therefore $\epsilon$ is vertical along $\gamma_1$,
and 
$\epsilon = cP_1$ for some constant $c$. However, we have just seen
that the pole solution $P_1$ has a nonzero horizontal component after
traversing the neck $c_1$. Thus the same argument shows $c=0$. 

If $f(\Sigma)$ is a cylinder then $P_1 = -P_2$ and the solution
$\epsilon = cP_1$ 
persists. Furthermore, $\epsilon$ remains exactly parallel to the
tangent vector as it traverses the radius $1/2$ circular
cross-sections of the cylinder, so it is tangent to $f(\Sigma^+)$. 
\hfill $\square$
\vskip .1in

We now consider a general coplanar $k$--unduloid with pole solutions
$P_1, \dots, P_k$. 

%In the case $k \geq 3$, one can still identify the polygon determined 
%by $P_1, \dots, P_k$. 

\begin {prop} \label {poles-prop}
The pole solutions $P_1, \dots, P_k$ are the vertices of the polygonal 
disc used in \cite{GKS2} to classify coplanar $k$-unduloids. 
\end {prop} 

\begin {rmk}
Within this proof, and for the remainder of the paper, we say a function
$u \simeq 0$ on an end $E_j$ if $u$ and its derivatives decay exponentially
on the end $E_j$. Similarly, a vector field $\epsilon \simeq 0$ on
an end $E_j$ if each of its components and their derivatives decay 
exponentailly.
\end {rmk} 

{\bf Proof}: Here, and later in Section \ref{neumann-sec}, we
truncate the symmetry curves $\gamma_j$ at approximate necks of the 
ends $E_j$ and $E_{j+1}$. By Lemma \ref{delaunay}, the 
length of $\gamma_j$ between successive necks is $\pi$, so each pole 
solution rotates through an angle of $2\pi$ from neck to neck. Thus the
value of the pole solutions is independent of which truncation of the
symmetry curves $\gamma_j$ we choose. We can compute the
distance in $S^2$ between $P_j$ and $P_{j+1}$ by traversing a neck curve 
$c_j$ of $E_j$, connecting $\gamma_j$ and $\gamma_{j+1}$. Exact unduloid
necks with the orientation indicated in Figure \ref{configuration} have 
conormal pointing in the axis $a_j$ direction, so along $c_j$ every
solution $\epsilon$ to equation (\ref{compass-eqn}) satisfies 
$$d\epsilon (c_j') = 2\epsilon \times df (J_0(c_j')) \simeq 
2\epsilon \times a_j.$$
This implies that (up to exponentially decaying terms, which are
negligible) each unit $\epsilon$ rotates with angular speed $2$ about
the $a_j$ axis as it traverses $c_j$. 
The total length of $c_j$ is $n_j/2$, so the total rotation
angle along $c_j$ is $n_j$. 
Choose positively oriented frames $\{ a_j, b_j, e_3\}$
for each end $E_j$,
as indicated in Figure \ref{configuration}. Then as we
traverse $c_j$ the pole solution $P_j$ rotates
 in a great circle of $S^2$, clockwise
in the plane spanned by $b_j$ and $e_3$, and we deduce that the
distance from $P_j$ to $P_{j+1}$ is $n_j$.  Thus the edge lengths of
the polygonal loop are exactly the necksizes $n_1,\dots,n_k$ of
$f(\Sigma)$. The polygonal disc used in \cite{GKS2} to classify 
coplanar $k$--unduloids also satisfies that the distance between successive 
vertices is $n_1, \dots, n_k$, so (after normalizing to fix $P_1$ and 
a frame at $P_1$ using a rotation) the two sets of vertices coincide.
\hfill $\square$
%Do we need to add anything to this? jesse

%Even in the Delaunay case ($k=2$) the $\epsilon$-fields contain
%useful information. 

There is an interesting consequence and generalization of the fact 
that the period of any unduloid is $\pi$.
Consider a
coplanar $k$--unduloid and let $L_j$ be the length 
of the curve $\gamma_j$ obtained by truncating at the (asymptotically
exact) necks $c_{j-1}$ and $c_j$. By Lemma \ref{delaunay}, the
length mod $\pi$ of these curves has a well-defined limit as the
truncations approach infinity. We call this limit $L_j^\infty$. 

\begin {prop} \label {2L_j}
Let $\alpha_j$ be the interior angle at the vertex $P_j$ of the
spherical polygon associated to 
$f(\Sigma)$, and let $\beta_j$ be the angle between the asymptotic axes
$a_{j-1}$ and $a_j$ (see Figure \ref{configuration}). Then 
$$2L_j^\infty = \pi + \alpha_j + \beta_j \mbox{ mod }2\pi.$$
\end {prop}

\begin {rmk} This result is equivalent to the relation found 
  (Proposition $7$ of \cite{GKS0}) for the twist angle of the conjugate
  cousin minimal surface around each of its boundary Hopf circles. 
\end {rmk}

\begin {figure}[h]
\begin {center}
\includegraphics[width=2.5in]{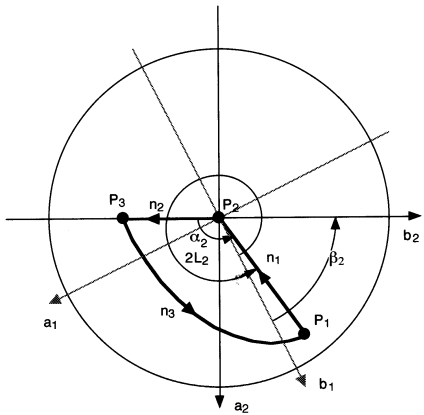}
\caption{The top view of the pole solutions just before traversing the
second neck}
\label {angles}
\end {center}
\end {figure}

{\bf Proof}: One can see from equation (\ref{compass-eqn}) that after
traversing $\gamma_j$, the horizontal components of the arc from
$P_{j-1}$ to $P_j$  have 
rotated through an angle $2L_j$. As indicated by the angle relations
illustrated in Figure \ref{angles} (for $j=2$ on a triunduloid), this
must be asymptotically equal (up to multiples of $2\pi$) to $\pi +
\alpha_j +\beta_j$. 
\hfill $\square$

\section {The proof of the main theorem} \label {nondegen-sec}
 
We prove Theorem \ref{main-thm} in this section.
The proof uses two features of the Alexandrov symmetry satisfied by a
coplanar $k$--unduloid $f:\Sigma \rightarrow \R^3$. First, the
reflection symmetry lets us decompose any Jacobi field $u$ into the
sum of an even part $u_+$ and an odd part $u_-$. We call an even field
Neumann 
because its restriction to $\Sigma^+$ satisfies 
$$\calL_f(u_+) = 0, \qquad \left. \frac{\del u_+}{\del \eta} \right|_{\del
  \Sigma^+} = 0,$$
where $\eta$ is the (outer) conormal to $\del \Sigma^+$. Similarly, we
  call an odd field Dirichlet since it vanishes on $\del
  \Sigma^+$. Second, the graphical nature of $f(\Sigma^+)$ implies
  that $v := -\langle \nu, e_3 \rangle$ is a positive Dirichlet Jacobi
  field on $f(\Sigma^+)$. Using $v$ as a comparison, we show in
  Section \ref{dirichlet-sec}
  that $0$ is the only $L^2$ Dirichlet Jacobi field. This analysis so
  far carries through for coplanar CMC surfaces of any genus.

In order to analyze the Neumann Jacobi fields in Section
\ref{neumann-sec}, we use the conjugate 
variation field $\tilde \epsilon$ constructed in Section
\ref{conj-var-sect}. This requires $\Sigma^+$ to be simply connected,
that is, $\Sigma$ must have genus zero. 

\subsection {Dirichlet Jacobi fields} \label {dirichlet-sec}

The proof we give of Proposition \ref{dirichlet} immediately below uses the 
maximum principle and is analogous to the standard proof that the first 
eigenvalue of $\Delta$ on a bounded domain $\Omega$ is simple. In Appendix 
\ref{2nd-dirichlet}
we prove a stronger version of Proposition \ref{dirichlet} using an integral 
version of the same maximum principle argument. 
Both proofs compare 
$u$ to the vertical translation field $v := -\langle \nu, e_3 \rangle
= -\nu_3$. Notice that $v>0$ 
on $\Sigma^+$ and $v = 0$ on $\del \Sigma^+$. 

To apply our maximum principle arguments comparing $u$ to $v$, we need
to know 
$$v_\eta := \frac{\del v}{\del \eta} \leq -\delta < 0$$
on $\del \Sigma^+$. (We continue our convention that $\eta$ is
the outer conormal, which in this case is $-e_3$ along $\del \Sigma^+$.)
One can quickly deduce this inequality for some positive $\delta$,
because it is true near the 
ends (with $\delta = 1$) and since on any compact subset of $\del
\Sigma^+$ the Hopf boundary point lemma gives a (noncomputable) value
for $\delta$. The following lemma shows that we may take $\delta =1$
along all of $\del \Sigma^+$. We include this lemma, which is a
reinterpretation of height and gradient estimates carried out in 
\cite{KKS, KK}, for its geometric consequences. 

\begin {lemma} \label {bndry-curv}
Let $f(\Sigma)$ be an Alexandrov symmetric (see Section \ref{modulisub})
CMC surface with finite
topology which is not a sphere. The boundary $\del f(\Sigma^+)$ is a
union of principal curves on $f(\Sigma)$ with principal curvature
$\kappa_1 < 1$. In particular,
the symmetry curves do not contain umbilics, and $\kappa_2 = \langle
\nu, e_3 \rangle_\eta =  -v_\eta > 1$.   
\end{lemma}

{\bf Proof}: Because $\del f(\Sigma^+)$ is the fixed point set of a
reflection symmetry for $f(\Sigma)$, it is a union of principal
curves. 

By the CMC equation, we have 
$$\Delta_f (z) = 2\nu_3,$$
where $z$ is the restriction of the vertical coordinate to the surface
$f(\Sigma^+)$. Also, because the components of the normal $\nu$ satisfy
the Jacobi equation, we have
$$\Delta_f(\nu_3) = -|A|^2 \nu_3 \geq -2 \nu_3,$$
Here we have used that $|A|^2 \geq 2$ and $\nu_3 < 0$. Thus
we have
$$\Delta_f(z + \nu_3) = (2-|A|^2) \nu_3 \geq 0,$$
and so $z + \nu_3$ is a subharmonic function on $\Sigma^+$. On $\del
\Sigma^+$, each 
function vanishes, so $z + \nu_3 = 0$. By explicit computation, $z +
\nu_3 \leq 0$ on the unduloid ends of $\bar \Sigma^+$. Thus in (the interior
of) $\Sigma^+$,  
$$ z + \nu_3 < 0$$
by the strong maximum principle. (Equality can only hold when $f$
parameterizes a unit hemisphere.) 

By the Hopf boundary point lemma, 
$$0 < \frac{\del}{\del \eta}(z + \nu_3) = -1 + \frac{\del
  \nu_3}{\del \eta} = -1 + \frac{\del}{\del \eta} \langle \nu, e_3
  \rangle.$$
We can rearrange this to obtain the curvature perpendicular to the
  boundary 
\begin {equation} \label {deriv-nu}
\kappa_2 = \frac{\del}{\del \eta} \langle \nu, e_3 \rangle  = -v_\eta >
     1, \end 
     {equation}
and so the principal curvature along the boundary is 
$$\kappa_1 = 2-\kappa_2 < 1.$$
\hfill $\square$

\begin {prop} \label{dirichlet}
Let $f(\Sigma)$ be a noncompact Alexandrov symmetric CMC surface 
of finite genus and with a finite number of ends. 
Then the only $L^2$ Dirichlet Jacobi field is the zero function. 
\end {prop}

{\bf Proof}: After possibly replacing $u$ with $-u$, we can 
assume $u>0$ somewhere. Now let $\mu > 0$ be a positive parameter. That $u 
\in L^2$ implies (see \cite {KMP}) that 
$u$ and its derivatives decay exponentially. Combining this
exponential decay with inequality (\ref{deriv-nu}), we see that for $\mu$
sufficiently large 
$$\mu v >  u$$
everywhere in the interior of $\Sigma^+$, with equality on
$\del\Sigma^+$. We define
$$\mu^* = \inf \{\mu > 0 \mbox{ }|\mbox{ } \mu v(p) >  u(p)
\, , p \in \Sigma^+ \}.$$
There is some finite $q$ which is a critical point of
$\mu^* v - u$ with critical value $0$. The point $q$ lies in either
the interior or the boundary of $\Sigma^+$. In both cases
$$u(q) = \mu^* v(q), \qquad \nabla u (q) = \mu^* \nabla v(q),$$
and $u \leq \mu^*v$ on $\Sigma^+$. In either case, the
strong maximum principle (the Hopf boundary point lemma if $q \in \del
\Sigma^+$) implies $u \equiv \mu^* v$. Because $u \in L^2$, this implies
$\mu^* = 0$ and thus $u \equiv 0$. 
\hfill $\square$

\subsection {Neumann Jacobi fields}
\label {neumann-sec} 

Given a Jacobi field $u$ on the coplanar $k$--unduloid
$f(\Sigma)$, the conjugate field $\tilde
\epsilon$ defined by equation (\ref{conj-variation}) yields a
conjugate Jacobi field $\tilde u := \langle \tilde \epsilon, \tilde
\nu \rangle$ on the surfaces $\tilde f(\Sigma^+)$ and
$f(\Sigma^+)$. By the correspondence $\tilde \epsilon = \tilde f
\epsilon$ relating solutions of equations (\ref{conj-variation}) and 
(\ref{revised-conj-variation}), we see 
$$\tilde u = \langle \tilde \epsilon, \tilde \nu \rangle = \langle
\tilde f \epsilon, \tilde f \nu \rangle = \langle \epsilon, \nu
\rangle.$$
By definition, $\tilde u$ is a Jacobi field on the conjugate cousin 
$\tilde f (\bar \Sigma^+)$ to the top half of $f(\Sigma)$. By 
Lemma \ref{same-Jacobi-op}, 
it is also a Jacobi field on $f(\bar \Sigma^+)$. 

Let $\calV$ denote the space of $L^2$ Jacobi fields on $f(\Sigma)$. 
Our plan is to convert even (Neumann) Jacobi fields $u \in \calV$ into $L^2$
Dirichlet Jacobi fields $\tilde u$, use Proposition \ref{dirichlet}
to deduce $\tilde u \equiv 0$, and use this to show $u \equiv 0$. In
order to carry out this procedure, $u$ must satisfy a finite number of
linear conditions, which is why Theorem \ref{main-thm} only bounds
the dimension of $\calV$, rather than asserting $\calV = \{ 0 \}$.

The solution $\epsilon$ to equation (\ref{revised-conj-variation}) 
is uniquely determined on $f(\bar \Sigma^+)$ once we choose an 
initial value $\epsilon(p)$ at some point $p$. We will choose 
$\epsilon =0$ at an endpoint of the truncated symmetry curve 
$\gamma_j$ discussed in the proof of Proposition \ref{poles-prop}. 
Because $u \simeq 0$ on all ends, solutions to equation 
(\ref{revised-conj-variation}) approach solutions to equation
(\ref{compass-eqn}) on the unduloid asymptote of each end $E_j$. Thus, 
provided $\gamma_j$ is sufficiently long, choosing $\epsilon = 0$ 
at an endpoint of the $\gamma_j$ on $E_j$ forces $\epsilon \simeq 0$ on 
$E_j$. We refer to this choice of initial value for $\epsilon$ as setting
$\epsilon = 0$ on $E_j$. 

By \cite{GKS2},
$f(\Sigma)$ has at least two non-cylindrical ends, one of which we label
$E_k$ (see Figure \ref{configuration}). From Lemma
\ref{delaunay} in Section \ref{modulisub}, a necessary condition for
attaining zero Dirichlet data on the end $E_k$ is that $\epsilon$ must
converge to $0$, and so we specify a unique conjugate field $\epsilon$
associated to $u$ by setting $\epsilon = 0$ on this non-cylindrical
end $E_k$. Starting at $E_k$, we compute how $\epsilon$ changes along
the contours $\gamma_j$, and along the ends $E_j$. 

By Lemma \ref{bndry-curv}, the $\gamma_j$ are principal curves, with
curvature $\kappa_1 < 1$ and constant conormal $-e_3$. We have seen by
Proposition \ref{dirichlet} that $u \in \calV$ is even. Thus we have 
\begin{eqnarray*}
d\epsilon (\gamma_j')& = & 2\epsilon \times df \circ J_0(\gamma_j') + df
\circ J_1 (\gamma_j') + d(u\nu)\circ J_0 (\gamma_j') \\
& = &-2\epsilon \times e_3 + u(\kappa_2-\kappa_1)f_{\eta} +
u_{\eta}\nu + u\nu_{\eta} \\
& = & -2\epsilon\times e_3 + u(\kappa_1-\kappa_2 +\kappa_2)e_3 + u_{\eta}\nu
\\
& = &  -2\epsilon\times e_3 + u\kappa_1e_3
\end{eqnarray*}
along $\gamma_j$. The geometric
interpretation of this equation is that the horizontal part of $\epsilon$
rotates about $e_3$, counterclockwise with speed $2$, and the 
vertical part of $\epsilon$ changes at a rate of $u\kappa_1$. 
Now set
\begin{equation}\label{heights}
h_j(u) := \int_{\gamma_j} d(\langle \epsilon, e_3 \rangle ) =
\int_{\gamma_j} u\kappa_1 ds,
\end{equation}
where $s$ is the arc-length parameter along $\gamma_j$. These {\em
  heights} 
$h_j(u)$ measure the change in the vertical components of $\epsilon$
as one traverses $\gamma_j$. They play a key role in our analysis. 

The integration defining the heights $h_j(u)$ associates a real number
to each symmetry curve $\gamma_j$. We encode this by 
defining the linear transformation $T:\calV \rightarrow \R^k$ by 
\begin {equation} \label {map-T}
T(u) = (h_1(u), \dots, h_k(u)).\end {equation}
\begin {prop} \label {zero-h-prop}
Let $f(\Sigma)$ be a coplanar $k$--unduloid, and let $\calV$ be the
space of $L^2$ Jacobi fields on $f(\Sigma)$. Then the linear
transformation $T:\calV \rightarrow \R^k$ defined by 
expression (\ref{map-T}) is injective. In particular, the dimension of
$\calV$ is at most $k$.  
\end {prop}

{\bf Proof}: We prove this proposition in two steps. First, show
that $T(u) = 0$ implies the conjugate Jacobi field $\tilde u$, which
is uniquely defined by our choice that $\epsilon = 0$ on the
non-cylindrical end $E_k$, must be identically zero. The second
step is to show that whenever $\tilde u \equiv 0$ then $u \equiv 0$. 

As we traverse  $\gamma_1$ from the end $E_k$ to the end $E_1$  only the
vertical part of $\epsilon$ changes, and the total change in this
component is $h_1(u) = 0$. Thus $\epsilon(p)$ converges exponentially to
$0$ on $\gamma_1$ as $p$ approaches infinity on the end $E_1$. Since
$\epsilon$ also converges to a homogeneous solution on $E_1$, we see
that $\epsilon$ converges to $0$ on the entire end $E_1$. Repeat this
argument successively, traversing $\gamma_j$ from $E_{j-1}$ to
$E_j$, using the hypothesis that each $h_j(u) = 0$. We deduce that
$\epsilon$ converges to $0$ exponentially along each end and that it
remains vertical along each $\gamma_j$. Thus $\tilde u = \langle
\epsilon, \nu \rangle$ decays exponentially to zero along each end and
is a Dirichlet field, because $\epsilon$ is vertical and $\nu$ is
horizontal along each $\gamma_j$. Therefore, after extending $\tilde
u$ to all of $f(\Sigma)$ by odd reflection, Proposition
\ref{dirichlet} implies $\tilde u \equiv 0$. 

We proceed to the second step, which we set  aside as a lemma. 
\begin {lemma} \label{lemma:u=0}
If the conjugate Jacobi field $\tilde u$ is identically zero, then so
is $u$. \end {lemma}

{\bf Proof}: We assume $\tilde u = \langle \tilde
\epsilon, \tilde \nu \rangle \equiv 0$, that is, the vector field $\tilde
\epsilon$ is tangent to $\tilde f (\Sigma^+)$. We pull
$\tilde \epsilon$ back to $\Sigma^+$ and denote its 
flow by $X_{\tilde\epsilon}(t)$. For small values of $t$, this is a
diffeomorphism $X_{\tilde\epsilon}(t): \Sigma^+ \rightarrow \Sigma^+$,
because $\epsilon$ is parallel to the conormal, and so $\tilde
\epsilon$ is tangent along $\del \tilde f (\Sigma^+)$.
Now define the one-parameter family of immersions 
$$\tilde f(t) = \tilde f \circ X_{\tilde\epsilon}(t): \Sigma^+
\rightarrow S^3.$$
This provides a family of reparameterizations of the minimal surface
$\tilde f(\Sigma^+) \subset S^3$. 

We produce a family of CMC surfaces $f(t)$ in
$\R^3$ by taking the conjugate cousin of this family of
reparameterization of $\tilde f (\Sigma^+)$. Rearrange the conjugate
family equation (\ref{conj-family}) to read
\begin {equation} \label {conj-cousin2}
df (t) = - \tilde f(t) ^{-1} d\tilde f(t) \circ J(t). \end {equation}
Using the inhomogeneous equation (\ref{conj-variation}) for $\epsilon$ 
and $\tilde f(t) = \tilde f + t \tilde \epsilon + O(t^2) = \tilde f (1 +
t \epsilon + O(t^2))$,
expand equation (\ref{conj-cousin2}) in powers of
$t$. One recovers $d(u\nu)$ as the $O(t)$ term in the expansion of
$df(t)$: 
\begin{eqnarray*}
df (t) & = & -(\tilde f (1 + t\epsilon))^{-1} [(d\tilde f) (1 +
  t\epsilon) + t\tilde f d\epsilon] \circ (J_0 + tJ_1) + O(t^2) \\ 
& = & -(1 - t\epsilon) \tilde f ^{-1} [d\tilde f \circ J_0 + t (
  (d\tilde f \circ J_0) \epsilon + d\tilde f \circ J_1 + \tilde f
  d\epsilon \circ J_0)] + O(t^2) \\ 
& = & - \tilde f^{-1} d\tilde f \circ J_0 + t [ \epsilon \tilde f^{-1}
  d\tilde f \circ J_0 - \tilde f^{-1} (d\tilde f \circ J_0) \epsilon -
  \tilde f^{-1} d\tilde f \circ J_1 - d\epsilon \circ J_0] + O(t^2)
  \\ 
& = & df + t[-\epsilon df + df \epsilon - df \circ J_0 \circ J_1 - (-
  \epsilon df + df \epsilon - d(u\nu) + df \circ J_1 \circ J_0)] +
  O(t^2) \\ 
& = & df + td(u\nu) + O(t^2). \end {eqnarray*}
We used the facts that $J_0^2 = -I$ and $J_0 \circ J_1 = -J_1 \circ
  J_0$ in the last steps. 

Integrate the one-form $df(t) = df + td(u\nu) + O(t^2)$ to
recover the immersion $f(t)$. In this integration we are free to
choose the value of $f(t)$ at a basepoint $p \in \Sigma^+$, and choose 
$f(t)(p) = f(p) + tu\nu(p)$. 
Then for any compact set $K \subset \Sigma^+$ and $q \in K$, we have 
\begin {eqnarray*}
f(t)(q) & = & f(p) + tu\nu(p) + \int_p^q df(t) \\
& = & f(p) + tu\nu(p) + \int_p^q d(f + tu\nu) + O(t^2)) = 
f(q) + tu\nu(q) + O(t^2). \end {eqnarray*}
However, this one-parameter family $f(t)$ is a conjugate cousin family 
for the fixed surface $\tilde f(\Sigma^+)$, so by Theorem 1.1 of \cite{GKS1},
the surfaces $f(t)$ can only vary by a family of
translations. Taking the derivative at $t=0$, this implies $u$ is the
normal part of an $\R^3$ translation, which implies $u \not \in
L^2$. Thus $\tilde u \equiv 0$ implies $u \equiv 0$, completing the
proof that $T$ is injective.
\hfill $\square$ 

\begin {prop} \label{imag-T}
Suppose $f(\Sigma)$ is a coplanar $k$--unduloid. Let $u \in \calV$,
and let $P_1,
\dots, P_k$ be the pole solutions to the homogeneous equation
(\ref{compass-eqn}) associated to the symmetry curves $\gamma_1,
\dots, \gamma_k$. Then for the constants $h_j := h_j(u)$, we have the
linear relation 
\begin {equation} \label {lin-dep} \sum_{j=1}^k h_j P_j \equiv 0 
\end {equation}
on $f(\Sigma^+)$. Thus, if the vertices of the classifying polygon 
for $f(\Sigma)$ span an $l$-dimensional subspace of $\R^3$, then $\calV$
is at most $(k-l)$--dimensional.
\end {prop}

{\bf Proof}: Let $\epsilon$ be the conjugate variation field which solves
equation (\ref{revised-conj-variation}) for the given $u\in \calV$, with
$\epsilon = 0$ on the end $E_k$. Traversing $\gamma_1$ from $E_k$ to
$E_1$, as in the previous proposition, we conclude that $\epsilon$
converges exponentially to the homogeneous solution $h_1 P_1$ on the
end $E_1$. Thus $\epsilon_1 = \epsilon - h_1 P_1$ solves equation
(\ref{revised-conj-variation}) with inital value $0$ on $E_1$, and
evolves along $\gamma_2$ with a vertical change of $h_2$. Thus
$\epsilon_1$ converges to the homogeneous solution $h_2 P_2$ along the
end $E_2$, so $\epsilon$ converges to $h_1 P_1 + h_2 P_2$ along this
end. Continuing this reasoning and traversing the remaining $\gamma_j$
in order, one returns to the end $E_k$, with $\epsilon$ converging to
the homogeneous solution $h_1 P_1 + \cdots + h_k P_k$. Since
$\epsilon$ is well-defined, this sum must be the initial asymptotic
homogeneous solution $0$. This shows the linear dependence
(\ref{lin-dep}). 

Evaluating the pole solutions at a point $q \in \Sigma^+$ yields
vertices for a representative classifying polygon for $f(\Sigma)$. The
linear relation (\ref{lin-dep}) implies that $(h_1, \dots, h_k)$
solves a 
homogeneous system of rank $l = \dim\Span \{ P_1(q), \dots, P_k(q) \}
\leq 3$. Since the solution space of this system 
is $(k-l)$--dimensional, and the linear transformation $T$ defined by
equation (\ref{map-T}) is injective, we conclude that $\dim
\calV \leq k - l$. 
\hfill $\square$

Using the fact that the vertices of the classifying polygon of a
coplanar $k$--unduloid span a two- or three-dimensional subspace of
$\R^3$ \cite{GKS2}, this completes the proof of Theorem \ref{main-thm},
and, as explained in the introduction, Corollary \ref{nondegen-thm}. 

\section {Extensions and open questions} \label {questions}

One can sharpen the proofs of Theorem \ref{main-thm} and Corollary
\ref{nondegen-thm} to show that
triunduloids with a cylindrical end are also nondegenerate. The theorem
below includes these triunduloids 
as a special case, and applies to a more general class of
$k$--unduloids. By Theorem $1.5$ of \cite {GKS2}, a coplanar
$k$--unduloid has at least two non-cylindrical ends. 
\begin {thm} \label{cyl-dim-bound}
Let $f(\Sigma)$ be a coplanar $k$--unduloid. If $f(\Sigma)$ has $d$
non-cylindrical ends and the vertices of the classifying polygon span
an $l$--dimensional subspace of $\R^3$, then the space $\calV$ of $L^2$
Jacobi fields has dimension at most $d-l$. In particular, if
$f(\Sigma)$ has exactly two 
non-cylindrical ends, or three non-cylindrical ends and classifying
polygon with vertices spanning $\R^3$, then it is nondegenerate. 
\end {thm}

{\bf Proof}: By Proposition
\ref{dirichlet}, any $u \in \calV$ is even, so we proceed as
in Section \ref{neumann-sec}. The key idea in the proof is
the observation (see 
Lemma \ref{delaunay}) that if $E_j$ is a cylindrical end, then
the pole solutions $P_j$ and $P_{j+1}$ are opposites, and are
asymptotically tangent along $E_j$. In other words, given $u \in \calV$
and a corresponding conjugate variation field $\epsilon$, if
$\epsilon$ is vertical along $\gamma_j$ then it is asymptotically
tangent on $E_j$ and continues to be vertical on
$\gamma_{j+1}$. Therefore, 
the conjugate Jacobi field $\tilde u = \langle \epsilon, \nu \rangle$ 
vanishes on $\gamma_j \cup \gamma_{j+1}$ and decays along $E_j$. More
generally, if $(E_r, \dots, E_{s-1})$ is a string of adjacent cylindrical
ends and $\epsilon$ is vertical on $\gamma_r$, then it is vertical on
all the symmetry curves $\gamma_r \cup \dots \cup \gamma_s$, implying 
$\tilde u$ vanishes on these symmetry curves and decays on the ends
$E_r, \dots, E_{s-1}$. 

We now develop the combinatorial tools needed to complete the
proof. The distribution of non-cylindrical ends on $f(\Sigma)$ leads
to a partitioning of the cyclically ordered set of symmetry curves
$(\gamma_1, \dots, \gamma_k)$ and their corresponding pole solutions
$(P_1, \dots, P_k)$ into substrings. Our substrings have the form $C
:= (\gamma_r, \gamma_{r+1}, \dots, \gamma_s)$, where the ends $E_r,
E_{r+1}, \dots, E_{s-1}$ are cylindrical while $E_{r-1}$ and $E_s$ are
not. In other words, $\gamma_r \cup \dots \cup \gamma_s$ connects the
non-cylindrical end $E_{r-1}$ to the next non-cylindrical end $E_s$,
through adjacent cylindrical ends. Notice that the singleton $C =
(\gamma_j)$ is a substring if neither $E_{j-1}$ nor $E_j$ are
cylindrical ends. 
Because each substring corresponds to a path
joining one non-cylindrical end to the next non-cylindrical end in the
cyclic ordering, the
total number of elements of the partition equals the number of
non-cylindrical ends $d$ on $f(\Sigma)$. 

If $C = (\gamma_r, \dots, \gamma_s)$ is a substring then, by
the previous discussion, the corresponding pole solutions $(P_r,
\dots, P_s)$ are all
parallel; in fact, for $r\leq j \leq s$, we have $P_j =
(-1)^{j-r}P_r = (-1)^{s-j}P_s$. Moreover, if $\tilde u$ decays on
$E_{r-1}$ and if 
$$\sum_{j = r}^s h_j (u) P_j = (\sum_{j = r}^s (-1)^{s-j} h_j(u))P_s =
0,$$ 
then $\tilde u$ vanishes on $\gamma_r \cup \gamma_{r+1} \cup \dots
\cup \gamma_s$ and $\tilde u$ also decays on the ends $E_r, \dots,
E_s$. We now define the linear transformation $\hat T: \calV \rightarrow
\R^d$ by 
$$\hat T(u) := (\hat h_1(u), \dots, \hat h_d(u)) := (\sum_{j =
  r_1}^{s_1}(-1)^{s_1 -j} h_j(u), \dots, \sum_{j =
  r_d}^{s_d}(-1)^{s_d -j} h_j(u)),$$ 
where the $m^{th}$ string of the cyclic partition is $(\gamma_{r_m},
\dots, \gamma_{s_m})$. If $\hat T(u) = 0$, then each alternating sum
$\hat h_m(u)$ is zero, and so $\tilde u$ is an $L^2$ Dirichlet Jacobi
  field. Lemma \ref{lemma:u=0} then implies $u \equiv
0$. Therefore, $\hat T$ is injective. 

The linear relation (\ref{lin-dep}) now reads 
$$0 \equiv \sum_{m= 1}^k h_m P_m = \sum_{m = 1}^d (\sum_{j = r_m}^{s_m}
(-1)^{s_m - j} h_j)P_{s_m} = \sum_{m = 1}^d \hat h_m P_{s_m}.$$
As in the proof of Lemma \ref{lemma:u=0}, this linear system has
rank $l = \dim 
\{ \Span\{P_{s_1}, \dots, P_{s_d}\} \leq 3$, so the solution space is
$(d-l)$--dimensional. Since $\hat T:\calV \rightarrow \R^d$ is
injective, we deduce that $\dim \calV \leq d-l$. 

\hfill $\square$
\vskip .1in

One can realize the space of triunduloids as a three-ball very 
explicitly in the following way (see also \cite{GKS1}). Evaluating the 
ordered triple of pole solutions $P_1,P_2, P_3$ at a base point yields a 
unique spherical triangle associated to the triunduloid. If the 
necksizes satisfy the strict triangle inequalities, this 
triangle is either strictly contained in an open hemisphere (which 
we normalize to be the upper hemisphere), or it strictly contains
a closed hemisphere (which we normalize to be the lower hemisphere). 
The vertices of this triangle are the pole solutions evaluated at our 
base point. In the first case, one can parameterize the vertices of all such 
triangles, and their associated triunduloids, by the upper half of an
open three-ball. Similarly, in the second case one can parameterize
the vertices of all such triangles and the associated triunduloids by
the lower half of an open three-ball.
Under this pair of 
parameterizations, the equatorial disc which joins these two half-balls
corresponds
to the triunduloids which satisfy the weak spherical triangle 
inequalities. 
By Corollary \ref{nondegen-thm}, each triunduloid corresponding to a 
point in the upper and lower half-ball under this parameterization is
nondegenerate. 

%One can realize the space of triunduloids as a three-ball very explicitly in
%the following way (see also \cite{GKS1}). When the necksizes satisfy 
%the strict spherical 
%triangle inequalities, the pole solutions lie in an open hemisphere, which we 
%normalize to be the upper hemisphere. There are 
%two spherical domains, bounded by the same triangle, associated to this triple
%of points on $S^2$, one of which is contained in the upper hemisphere. 
%Varying the pole solutions within the upper hemisphere, one can parameterize
%the resulting space of spherical domains in the upper hemisphere by the 
%upper half of an open three-ball.
%Similarly, one can parameterize the spherical domains 
%which extend 
%into the lower hemisphere by the lower half. 

%Using Proposition
%\ref{cyl-dim-bound}, we can show: 
\begin{cor}
The nondegenerate triunduloids form a connected open subset in
the space of all triunduloids. 
\end{cor}

{\bf Proof}: %Observe that, by the Implicit Function Theorem, the set of
%nondegenerate CMC surfaces is open.
Observe that by the Implicit Function Theorem, the set of nondegenerate
triunduloids is open. To 
show connectedness, it suffices to find a
nondegenerate triunduloid satisfying the weak spherical triangle
inequalities, which lies in the closure of the two open half-balls 
described above. Any triunduloid with a cylindrical end is such a 
surface. 
\hfill $\square$
\vskip .1in

%\subsection {Open questions}

We conclude by mentioning several naturally related open problems
concerning Jacobi fields on CMC surfaces and the moduli space theory
of CMC surfaces. Theorems \ref{main-thm} and \ref{cyl-dim-bound} give
upper bounds for the dimension of the 
space $\calV$ of $L^2$ Jacobi fields on coplanar $k$--unduloids. Is
this bound sharp? In particular, up to scaling, there is at most
one nonzero $L^2$ Jacobi field on any triunduloid satisfying $n_1 + n_2 +
n_3 = 2\pi$ or $n_i + n_j = n_k$. Does this Jacobi field ever exist? 

Is it possible to extend our technique to a wider class of CMC
surfaces? For instance, there are many CMC surfaces which are not
Alexandrov symmetric but do have some symmetry ({\it e.g.} tetrahedral
symmetry). Can one use our methods to bound either the necksizes or
the dimension of $\calV$ on such
surfaces? Might the analysis of Section \ref{neumann-sec} also 
bound the dimension of $\calV$ on
Alexandrov-symmetric CMC surfaces with positive genus? 

It would be very interesting to produce an example of a 
degenerate CMC surface. The question of integrability of a Jacobi
field is also open. According to \cite {KMP}, any tempered
(sub-exponential growth) Jacobi field on a 
nondegenerate CMC surface is integrable, in the sense that it is the
velocity vector field of a one-parameter family of CMC surfaces. It
would be useful to decide whether tempered
Jacobi fields are always integrable in this sense. 

\appendix
\addappheadtotoc
\appendixpage

\section {The spinning sphere connection} \label{spin-sphere-appen}

Some of the material in this section is well known, but we include it for the 
convenience of the reader. We begin with a proof of Proposition 
\ref{spin-sphere-prop}:

{\em Let $f:\Omega \rightarrow \R^3$ be an immersion and 
consider the $SO(3)$-connection defined by spinning a sphere at
speed $2$, as described in Section \ref{spin-sphere-sec}. The $f(\Omega)$ 
has mean curvature $1$ if and only if this connection is flat.}

{\bf Proof}: We have already shown 
that the CMC condition implies 
the flatness; it remains to prove the reverse implication. The
assumption that the spinning sphere connection is flat is exactly the
hypothesis that equation (\ref{compass-eqn}) is integrable for
$\epsilon$ on any simply connected domain $\Omega$, for any choice of
initial vector $\epsilon(f(p))$. Using equation (\ref{compass-eqn}),
integrability implies  
$$0 = d(d\epsilon) = 2 [2(\epsilon \times df \circ J_0) \times df
  \circ J_0] + 2\epsilon \times (d(df \circ J_0)).$$
The second term is $2\epsilon
  \times (-\Delta_0 f) dx \wedge dy$. Expand the first term
  and then use the Jacobi identity: 
\begin {eqnarray*}
4 (\epsilon \times df \circ J_0) \times df \circ J_0  & = & 4
(\epsilon \times (f_y dx - f_x dy)) \times (f_y dx - f_x dy) \\ 
& = & 4 (-(\epsilon \times f_y) \times f_x + (\epsilon \times f_x)
\times f_y) dx \wedge dy = 4 \epsilon \times (f_x \times f_y) dx
\wedge dy .
\end {eqnarray*}
Now combine these two terms to obtain 
$$0 = d(d\epsilon) = 2\epsilon \times (2f_x \times f_y - \Delta_0 f)
dx \wedge dy.$$
Because $\epsilon$ can be chosen to have any value at a point,
we deduce that $f$ solves equation (\ref{cmc-eqn}). 
\hfill $\square$
\vskip .1in

%\begin {rmk}
%Amusingly, this even
%gives a nontrivial property of the round sphere -- the Earth, for
%example; in principle, if one were to keep track of the orientation of
%a spinning half-sized sphere -- for example, Mars -- this would give a
%kind of ``global positioning system''; perhaps NASA should look into
%this -- if ``Sisyphus'' Bush can't go to Mars, perhaps he could bring it
%back to roll around the Earth?!
%\end {rmk}

The solutions $\epsilon$ to the homogeneous system 
 (\ref{compass-eqn})
can also be expressed naturally in terms of
 the quaternion geometry of $S^3$ and the conjugate
surface equation for 
 $\tilde f:\Omega\rightarrow S^3.$ Following the ideas in the
 abstract sketch of the proof of Propsition \ref{conj-var-exist}, 
let 
$$q(t)=1+t\alpha+O(t^2)$$
be a smooth 
 curve of unit quaternions, passing through $1$ at time $t=0$, with  
$\alpha\in T_1 S^3 = \R^3$, a fixed imaginary quaternion. Consider
the family of left translations $q(t)\tilde f$ of the mapping $\tilde f$,
and note that since the translation isometry is on the left, each of
these surfaces satisfies the conjugate cousin equation,
$d(q(t)\tilde f) = (q(t)\tilde f) df\circ J_0.$
Therefore, the velocity $\tilde \epsilon = \alpha \tilde f$ of the
family at $t=0$ solves the homogeneous ($u \equiv 0$) version of equation
(\ref{conj-variation}), and 
\begin {equation} \label {left-right}
\epsilon := \tilde f^{-1} \tilde \epsilon = \tilde f ^{-1} \alpha
\tilde f \end {equation}
solves equation (\ref{compass-eqn}). (One can also check by direct
computation that $\epsilon = \tilde f^{-1} \alpha \tilde f$ solves equation
(\ref{compass-eqn}).) By varying $\alpha$ one obtains in this manner
the unique solution to each 
initial value problem for equation (\ref{compass-eqn}).

Continuing our interpretation of equation (\ref{left-right}), we see 
that an equivalent way to understand
the spinning-sphere
flat connection on $f(\Omega)$ is as
the pullback from $\tilde f(\Omega)$ to $f(\Omega)$ of a natural double
covering  $S^3\rightarrow SO(3)$, arising from quaternion conjugation:
for each imaginary quaternion $\alpha\in\R^3$ and each $q\in S^3$,
write
\begin{equation} \label{SO(3)-eqn}
R_{q}(\alpha) := q^{-1}\alpha q.
\end{equation}
We have seen that for fixed $\alpha$ the $\R^3$-valued field on $S^3$
defined by equation (\ref{SO(3)-eqn}) pulls back to a
solution of equation (\ref{compass-eqn}) on $f(\Omega)$. 
More generally,
for each $q\in S^3$  the linear map $R_q$
 is actually a
rotation (in $SO(3)$), and the flat connection on
$f(\Omega)$ is the pullback of this rotation field from $S^3$.

%Actually the involuted conjugation
%map $q\rightarrow R_{q^{-1}}$ is a double
%covering homomorphism from $S^3$ to $SO(3)$.  This is easy to see 
%directly from equation (\ref{SO(3)-eqn}):  the identity
%rotation $R_q=I$ arises if and only if the unit $q$ commutes with all
%quaternions, which is equivalent to $q\in\{1,-1\}.$  The homomorphism
%property then implies that $R_{q_1}=R_{q_2}$
%if and only if $R_{q_1(q_2)^{-1}}=I$, that is, $q_1,q_2$ are equal or
%opposite.

%One can check that $q\rightarrow R_q$ is onto
%$SO(3)$ by
%explicitly computing the rotation given by $R_q$.
%If we write  $q=\Exp(t\beta)$, where $\beta$ is a unit
%imaginary quaternion, then we claim that
%$R_q$ is a rotation
%with axis $\beta$, and that $R_q$ rotates
%an amount $-2t$ in the positive direction about the $\beta$- axis.
%Quaternion algebra  verifies these claims. 
%First verify that  $\beta$ is fixed by $R_q$:
%$$R_{\Exp(t\beta)}(\beta)=\Exp(-t\beta)\beta \Exp(t\beta)=\beta,$$
%because the three terms in the product commute.  Next consider an
%imaginary quaternion $\alpha$
%perpendicular to $\beta$:
%\begin {eqnarray*}
%R_{\Exp(t\beta)}(\alpha) & = & \Exp(-t\beta)(\alpha \Exp(t\beta))\\
%& = & \Exp(-t\beta)(\Exp(-t\beta)\alpha) = \Exp(-2t\beta)\alpha\\
%& = &  \cos(2t)\alpha - \sin(2t)(\beta\times\alpha).
%\end {eqnarray*}
%Since $\{\alpha,\beta\times\alpha,\beta\}$ is positively oriented,
%it follows that  $\alpha$ is rotated by an angle $-2t$ 
%about the $\beta$ axis, as claimed.

We conclude from this discussion that the rotation of the spinning
sphere 
$$R := R_{\tilde f} : \Omega \rightarrow SO(3)$$ 
is nothing more than the conjugate cousin $\tilde f$ followed by the
natural covering map $S^3 \rightarrow SO(3)$. 
Because $\tilde f$ is harmonic, so is the map $R$. (One can
verify this directly using (\ref{compass-eqn}) to compute
$$ R^{-1} \Delta_0 R = (R^{-1} R_x)^2 + (R^{-1} R_y)^2,$$
which is the equation for a harmonic map from $\Omega \subset
\R^2$ to $SO(3)$, see \cite{U}). Furthermore, the solution
$\epsilon$ to equation (\ref{left-right}) is $R(\alpha)$. 

%\section {The second appendix}

\section {Bounded Dirichlet Jacobi fields}
\label {2nd-dirichlet}

\begin {prop} 
Let $f(\Sigma)$ be an Alexandrov symmetric CMC surface (see Section
\ref{modulisub}) of finite genus and with a finite number of ends. 
Every bounded odd (Dirichlet) Jacobi field $u$ on $f(\Sigma)$ is a
constant multiple of the vertical translation field $v = -\langle \nu,
e_3 \rangle =-\nu_3$. In particular, if $u \in \calV$ then $u$ is an
even (Neumann) Jacobi field (unless $f(\Sigma)$ is a unit sphere). 
\end {prop}

{\bf Proof}: We initially assume $u \in \calV$, rather than the
more general hypothesis that $u$ is a bounded Dirichlet Jacobi
field. For this proof it is technically simpler to consider the entire
surface $f(\Sigma)$. Recall that both $u$ and $v$ are odd with
respect to reflection through the Alexandrov plane of symmetry, and by
inequality (\ref{deriv-nu}) $u/v$ is uniformly bounded on the
complement of the symmetry curves, which is $\{v \neq 0\}$. Also, both
$u$ and $v$ are real analytic functions which vanish on the symmetry
curves. These facts imply that $u/v$ extends to an even, real analytic
function on the entire surface $f(\Sigma)$. To verify analyticity on
$\{ v = 0 \}$, use conformal
curvature coordinates in which the $x$--axis is a symmetry curve; the
fact that $u$ and $v$ both vanish on the $x$--axis means we can write 
$$u(x,y) = y U(x,y), \qquad v(x,y) = y V(x,y), \qquad
\frac{u(x,y)}{v(x,y)} = \frac{yU(x,y)}{yV(x,y)} = 
\frac{U(x,y)}{V(x,y)}, $$
where $U$ and $V$ are also real analytic and $V \neq 0$ near the
$x$--axis by Lemma \ref{bndry-curv}. 

Continuing with the proof, assume that 
$u/v >0$ somewhere. Since $u/v$ is nonconstant, we can pick a
regular value $\delta > 0$ for $u/v$ 
with nonempty inverse image. The domain 
$$\Omega_\delta := \{ u/v > \delta \}$$
is bounded (because $u \in L^2$) and has smooth boundary in $\Sigma$. 
Since $(u/v)_\eta< 0$ pointwise along $\del \Omega_\delta$, 
\begin {equation} \label {neg-integrand}
\int_{\del \Omega_\delta}  v \frac{\del u}{\del \eta} - u \frac{\del
  v} {\del \eta} = \int_{\del \Omega_\delta} v^2 \frac{\del
  (u/v)}{\del \eta}  < 0.\end {equation}
However, we also have 
\begin {eqnarray} \label {int-by-parts1}
0 & = & \int_{\Omega_\delta} v\calL_f u - u \calL_f v = \int_{\Omega_\delta}
v\Delta_f u - u \Delta_f v \\ \nonumber
& = & \int_{\del \Omega_\delta} v \frac{\del u}{\del \eta} - u
\frac{\del v}{\del \eta} = \int_{\del \Omega_\delta} v^2 \frac{\del
  (u/v)}{\del \eta}. \end {eqnarray}
This last equation (\ref{int-by-parts1}) contradicts the previous
inequality (\ref{neg-integrand}), proving $u \equiv 0$. 

We now explain how to extend this argument to prove that any bounded
Dirichlet Jacobi field $u$ is a constant multiple of $v$. We assume
$u/v$ is nonconstant and positive somewhere, pick a regular
value $\delta>0$, and define the nonempty set $\Omega_\delta$ as
before. In this case the inequalty (\ref{neg-integrand}) still holds,
but we cannot immediately deduce equation (\ref{int-by-parts1})
because $\Omega_\delta$ may be unbounded. We
overcome this difficulty by appealing to the linear decomposition
lemma of \cite{KMP}, which implies that on each end $E_j$, we have
exponential convergence 
$$u \simeq \sum_{i=1}^3 a_{ij} \nu_i,$$
where $\nu_i$ are the components of the normal vector to the
asymptotic unduloid. (In the case when the end $E_j$ is cylidrical,
one must also include Jacobi fields arising from changing the
necksize, which are even.) Because $u$ is odd, we must have $u \simeq a_j
\nu_3 :=  a_{3j} \nu_3$ on the end $E_j$, and so $u/v$ converges
smoothly to a constant $-a_j$ on the end $E_j$. ({\it A priori}, these 
constants may differ from end to end.) 

Now we truncate the domain $\Omega_\delta$ by intersecting $f(\Sigma)$
with a sequence of balls, defining 
$$\Omega_{\delta, N} :=\{ p \in \Omega_\delta : |f(p)|
\leq N \} = \Omega_\delta \cap \bar B_N (0),$$ 
where $N = 1,2,3, \dots$. Then equation (\ref{int-by-parts1}) becomes 
\begin {equation} \label {int-by-parts2}
0 = \int_{\Omega_{\delta,N}} u \calL_f v - v \calL_f u = \int_{\del
  \Omega_\delta \cap B_N} u v_\eta - v u_\eta + \int_{
  \Omega_\delta \cap \del B_N} u v_\eta - v u_\eta. \end {equation}
But as soon as $N$ is large enough so that $\del \Omega_\delta \cap 
  B_N$ has positive length, inequality (\ref{neg-integrand}) implies
  the first term is negative, and in fact it is decreasing in
  $N$; also, the second terms converge uniformly to zero by our
  previous discussion of the asymptotics.  This contradiction shows
  $u$ is a constant multiple of $v$. 
\hfill $\square$

\section {More about gluing}

In this section we outline another gluing construction which produces CMC 
surfaces with no small necks. 

The end-to-end gluing construction (Theorem $1$ of \cite{R}) proceeds 
as follows. Suppose
$f_1(\Sigma_1)$ and $f_2(\Sigma_2)$ are two nondegenerate CMC surfaces
with ends 
$E_j \subset f_j(\Sigma_j)$, such that $E_1$ and $E_2$ are asymptotic
to congruent Delaunay unduloids which are not cylinders. 
We must also assume that $f_1$ belongs to a one-parameter
family of CMC surfaces which changes the necksize of $E_1$ to first
order. Under these assumptions,
one can truncate $f_1(\Sigma_1)$ and $f_2(\Sigma_2)$ at necks of $E_1$
and $E_2$ and, after perturbation, glue together the resulting
surfaces with boundary to obtain a new CMC surface. The resulting CMC
surface is nondegenerate 
and has asymptotics which are close to the asymptotics of the
remaining ends of $f_1$ and $f_2$. One particular instance of the
end-to-end gluing construction, doubling along an end, 
occurs when one glues $f(\Sigma)$ to a copy of itself after
truncating a particular end. 

By Corollary \ref{nondegen-thm}, one can use most triunduloids in
end-to-end gluing, and in many other gluing constructions. In particular,
let $f(\Sigma)$ be a triunduloid which is a regular point of the 
classifying map and has necksizes $n_1, n_2, n_3$ 
satisfying the strict spherical triangle inequalities. By Sard's theorem, 
except for
a set of measure zero, all triunduloids with necksizes satisfying the strict 
triangle inequalities are regular points of the classifying map. Any end of 
any such $f(\Sigma)$ satisfies all the hypotheses for end-to-end gluing, and 
so one can double such a triunduloid along any of its ends. 
This gluing construction yields examples of 
nondegenerate $k$--unduloids with $k > 3$ and no small necks (that is, no
short closed geodesics). In addition, one can use end-to-end
gluing to create nondegenerate CMC surfaces with any finite topology
and no small necks. 

\section {Comparison of the CMC and minimal cases}

We now compare our proof and Cos\'in and Ros' \cite{CR} proof of the 
analogous result for genus zero, coplanar, minimal $k$--noids. Because of 
the special properties of finite total curvature minimal surfaces, they 
are able to prove that 
all bounded Jacobi fields on $f$ are linear combinations of the components of 
the unit normal vector $\nu$. 

A sketch of their proof proceeds as follows. Let $\calW$ be the space
of bounded Jacobi fields on a genus zero, coplanar, minimal $k$--noid
$f : \Sigma \rightarrow \R^3$. As in the CMC case, $f$ is Alexandrov
symmetric, so one can decompose $\calW$ into its even (Neumann) and
odd (Dirichlet) parts.  First pull back the round metric on $S^2$ to
$\Sigma^+$ using the Gauss map. Because $f(\Sigma)$ is minimal, this
yields a metric conformal to the induced metric on $\Sigma^+$,
accomplishing two things: it compactifies $\Sigma^+$, identifying the
ends to points, and it transforms the Jacobi operator into $\Delta_1 +
2$, where $\Delta_1$ is the Laplacian in the round metric. The
uniqueness of Dirichlet fields (up to scaling) now follows from the
fact that $v = - \langle \nu, e_3 \rangle$ is positive on
$\Sigma^+$. Next, given any bounded Jacobi field $u$ on a finite total
curvature minimal surface, one can construct an associated branched
minimal surface $X(\Sigma)$ with planar ends which has the same Gauss
map as $f(\Sigma)$, and has $u$ as its support function (inner product
of the position vector and unit normal vector, that is, the Jacobi
field corresponding to the invariance of the minimal surface equation
under homothety).  The conjugate surface $\tilde X(\Sigma)$ is again
a minimal surface with planar ends, and its support function $\tilde
u$ again is a bounded Jacobi field. Since this conjugation from $u$
to $\tilde u$ interchanges Neumann and Dirichlet fields, Cos\'in and
Ros conclude that the Neumann fields on $f(\Sigma)$ are
(multiples of) the horizontal components of the unit normal.

One can also prove their result using our methods. To prove the uniqueness of 
Dirichlet Jacobi fields, up to scaling, one can slightly modify our proof in 
Appendix \ref{2nd-dirichlet}. The salient feature one must recall is that any 
bounded Jacobi field $u$ has a decomposition on each end $E$ as 
$$u = a_0 u_0 + a_1 u_1 + a_{-1} u_{-1} + O(r^{-2}),$$
where $r$ is the Euclidean distance
from the axis of the catenoid asymptote of $E$, $u_0 = O(1)$ arises
from translation along the asymptotic axis, and $u_{\pm 1} =
O(r^{-1})$ arise from translations perpendicular to the asymptotic
axis. We will make the normalization that $u_1$ corresponds to vertical 
translations. 
In particular, if $u \in \calW$ is Dirichlet, then $u = a_1
\langle \nu, e_3 \rangle + O(r^{-2})$, and so the boundary terms in
equation (\ref{int-by-parts2}) caused by spherical truncation approach
zero. Thus every bounded Dirichlet Jacobi field is a constant multiple
of $\langle \nu, e_3 \rangle$. 
One can then transform Neumann Jacobi fields to Dirichlet Jacobi fields 
using a construction analogous to our conjugate Jacobi field construction 
of Section \ref{conj-var-sect}. In this case, the conjugate 
variation field $\epsilon$ satisfies 
\begin {equation} \label {min-conj-var}
d\epsilon = df \circ J_1 + d(u\nu) \circ J_0.
\end {equation}
Now argue as in Section
\ref{neumann-sec}, using the heights $h_j(u)$ defined by equation
(\ref{heights}), which still measure the vertical change in $\epsilon$
evolving by equation (\ref{min-conj-var}) along $\gamma_j$. Because
equation (\ref{min-conj-var}) contains no rotation term and $d\epsilon
= O(r^{-2})$ on the ends, $\epsilon$ remains vertical along all the
symmetry curves and at infinity. Thus $\tilde u = \langle \epsilon,
\nu \rangle$ is a bounded Dirichlet Jacobi field, and we apply the
proof of Proposition \ref{zero-h-prop} to conclude $u = \langle \nu, b
\rangle$ for some $b \in \R^3$.


\begin {thebibliography}{99}

\bibitem [A]{A} A. D. Alexandrov. {\em A characteristic property of
  spheres}. Ann. Mat. Pura Appl. 58:303--315, 1962. 

\bibitem [CR] {CR} C. Cos\'in and A. Ros. {\em A Plateau problem at
  infinity for properly immersed minimal surfaces with finite total
  curvature}. Indiana Univ. Math. J. 50:847--879, 2001. 

\bibitem [G]{G} K. Gro\ss e-Brauckmann. {\em New surfaces of constant
  mean curvature}. Math. Z. 214:527--565, 1993. 

\bibitem [GKS0]{GKS0} K. Gro\ss e-Brauckmann, R. Kusner and
  J. Sullivan. {\em Classification of embedded constant mean curvature
  surfaces with genus zero and three ends}. Bonn SFB 256, Preprint
  536, 1997. 

\bibitem [GKS1] {GKS1} K. Gro\ss e-Braukmann, R. Kusner and J. Sullivan
  {\em Triunduloids: Embedded constant mean curvature surfaces with
  three ends and genus zero}. J. Reine Angew. Math. 564:35--61, 2003.

\bibitem [GKS2]{GKS2} K. Gro\ss e-Brauckmann, R. Kusner and
  J. Sullivan. {\em Coplanar constant mean curvature surfaces}. Preprint. 

\bibitem [Ho]{Ho} H. Hopf. {\em Differential Geometry in the
  Large}. Lecture Notes in Mathematics 1000. Springer Verlag. 1983.

\bibitem [Kap]{Kap} N. Kapouleas. {\em Complete constant mean curvature 
surfaces in 
Euclidean three-space}. Ann. Math. 2:239--330, 1990.

\bibitem [K]{K} H. Karcher. {\em The triply periodic minimal surfaces
  of Alan Schoen and their constant mean curvature
  companions}. manuscripta math. 64:291--357, 1989.

\bibitem [KK]{KK} N. Korevaar and R. Kusner. {\em The global structure
  of constant mean curvature surfaces}. Invent. Math. 114: 311--332, 1993

\bibitem [KKS]{KKS} N. Korevaar, R. Kusner and B. Solomon. 
{\em The Structure of Complete Embedded Surfaces with Constant Mean
Curvature}. J. Differential Geom. 30:465--503, 1989.

\bibitem [KMP]{KMP} R. Kusner, R. Mazzeo and D. Pollack. {\em The
  moduli space of complete embedded constant mean curvature
  surfaces}. Geom. Funct. Anal. 6:120--137, 1996.

\bibitem [L]{L} H. B. Lawson. {\em Complete minimal surfaces in
  $S^3$}. Ann. Math. 92:335--374, 1970. 

\bibitem [MP]{MP} R. Mazzeo and F. Pacard. {\em Constant mean curvature 
surfaces with
Delaunay ends}. Comm. Anal. Geom. 9:169--237, 2001. 

\bibitem [MPP1]{MPP1} R. Mazzeo, F. Pacard and D. Pollack. {\em Connected 
sums of 
constant mean curvature surfaces in Euclidean 3 space}. J. Reine Angew. Math. 
536:115--165, 2001. 

\bibitem [MPP2] {MPP2} R. Mazzeo, F. Pacard and D. Pollack. {\em The
  conformal theory of Alexandrov embedded constant mean curvature
  surfaces in $\R^3$}. To appear in Proceedings of the Clay Summer
  Institute on the Global Theory of Minimal Surfaces. arXiv:math.DG/0110099

%\bibitem [O]{O} R. Osserman. {\em A survey of minimal
%  surfaces}. Dover, 1986. 

\bibitem [R]{R} J. Ratzkin. {\em An end-to-end gluing construction for
  surfaces of constant mean curvature}. Ph. D. thesis,
  Univ. Washington, 2001. 

%\bibitem [Sc]{Sc} R. M. Schoen. {\em Uniqueness, symmetry and
%  embeddedness of minimal surfaces}. J. Differential
%  Geom. 18:791--809, 1983.

\bibitem [Sm]{Sm} N. Smale. {\em A bridge principle for minimal and 
constant mean 
curvature submanifolds of $\R^N$}. Inv. Math. 90:505--549, 1987.

\bibitem [U]{U} K. Uhlenbeck. {\em Harmonic maps into Lie groups
  (classical solutions to the Chiral problem)}. J. Differential
  Geom. 30:1--50, 1989. 

\end {thebibliography}
\end {document}